\documentclass[11pt]{amsart}
\input{xypic}
\usepackage{amsmath,amsfonts,amssymb,amsthm}
\usepackage{rotating}
\newtheorem{theorem}[equation]{Theorem}
\newtheorem{proposition}[equation]{Proposition}
\newtheorem{corollary}[equation]{Corollary}
\newtheorem{definition}[equation]{Definition}
\newtheorem{observation}[equation]{Observation}
\newtheorem{lemma}[equation]{Lemma}
\numberwithin{equation}{section}
\theoremstyle{remark}
\newtheorem{remark}[equation]{Remark}
\newtheorem{notation}[equation]{Notation}
\newtheorem{example}[equation]{Example}
\DeclareMathOperator{\GL}{\mathrm{GL}}
\DeclareMathOperator{\Ker}{\mathrm{Ker}}
\DeclareMathOperator{\Id}{\mathrm{Id}}
\DeclareMathOperator{\Res}{\mathrm{Res}}
\DeclareMathOperator{\Hom}{\mathrm{Hom}}
\DeclareMathOperator{\Tr}{\mathrm{Tr}}
\DeclareMathOperator{\Trace}{\mathrm{Trace}}
\DeclareMathOperator{\Gal}{\mathrm{Gal}}
\DeclareMathOperator{\Ind}{\mathrm{Ind}}
\DeclareMathOperator{\Ref}{\mathrm{Ref}}

\newcommand{\reg}{{\mathrm{reg}}}
\newcommand{\Gut}{{\Psi}}
\DeclareMathOperator{\dotequal}{{\hskip0.2em\dot{=}\hskip0.2em}}
\newcommand{\BA}{{\mathbb A}}

\newcommand{\BQ}{{\mathbb Q}}

\newcommand{\C}{{\mathbb C}}
\newcommand{\bZ}{{\mathbb Z}}
\newcommand{\BZ}{{\mathbb Z}}
\newcommand{\BN}{{\mathbb N}}
\newcommand{\CA}{{\mathcal A}}
\newcommand{\CH}{{\mathcal H}}
\newcommand{\CB}{{\mathcal B}}
\newcommand{\CN}{{\mathcal N}}

\newcommand{\CM}{{\mathcal M}}
\newcommand{\CF}{{\mathcal F}}

\newcommand{\CV}{{\mathcal V}}
\newcommand{\lr}{\longrightarrow}
\newcommand{\inv}{^{-1}}
\newcommand{\lexp}[2]{{}^{#1}#2}

\newcommand{\cf}{{\it cf}}
\newcommand{\genby}[1]{\langle#1\rangle}
\newcommand{\Nv}{\CN_{\genby v}}
\newcommand{\e}{\varepsilon}
\newcommand{\eps}{\varepsilon}
\newcommand{\ogr}[1]{\mathop\otimes\limits^{\mathrm{gr}}#1}

\renewcommand{\le}{\hspace{0.1em}\mathop{\leqslant}\nolimits\hspace{0.1em}}
\renewcommand{\ge}{\hspace{0.1em}\mathop{\geqslant}\nolimits\hspace{0.1em}}

\begin{document}
\author{C.~Bonnaf\'e, G.I.~Lehrer and J.~Michel}
\title[Reflection cosets]{Twisted invariant theory for reflection groups}
\date{11th February, 2006}
\maketitle
\centerline {\em Dedicated to George Lusztig on his 60th birthday.}
\tableofcontents

\section{Introduction}

Let $G$ be  a finite group generated by (pseudo)reflections in  a vector space
of  dimension $r<\infty$  over the  algebraically closed  field $K$  of
characteristic  zero.  The purpose  of  this  work  is to  discuss  the
``twisted case'' of various phenomena associated with the structure and
invariant theory of $G$. That  is, we take an element $\gamma\in\GL(V)$
which normalises  $G$, and consider how  it acts on the  invariants and
covariants (for various representations) of  $G$, and properties of its
eigenspaces. In particular  we study generalisations of  the results of
\cite{LS1} and  \cite{LM}, and the action  of $\langle G,\gamma\rangle$
on  covariants.  Our basic  method  is  a  variation  on the  theme  of
\cite{L},  which  enables  us  to  relate  various  sets  of  constants
associated to $G$ with the  corresponding ones for parabolic subgroups.
In many  cases, we  identify the relevant  constants as  eigenvalues of
certain transformations.

Specifically, to  any finite dimensional  $\genby{G,\gamma}$-module $M$
we  associate a  (multi)set  of $m=\dim  M$ constants  $\eps_\iota(M)$,
which  depend  only  on  the coset  $G\gamma$.  Various  algebraic  and
geometric  properties  of  $G\gamma$  may then  be  expressed  in  terms
of  the   constants  $\eps_\iota$.  For  example   the  condition  that
$\zeta\in  K^\times$  be  regular  for  $G\gamma$  may  be
expressed in  terms of  the $\eps_\iota(V)$ and  $\eps_\iota(V^*)$ (see
\S\ref{regularity} below).  A theme of this  work is that if  $G'$ is a
parabolic subgroup  of $G$  which is normalised  by $\gamma$,  then the
constants $\eps_\iota(M)$  are the  same whether $M$  is regarded  as a
module  for $\genby{G,\gamma}$  or  $\genby{G',\gamma}$ 
(see (\ref{eqlists}) below.  This idea  is
behind many of the results and their proofs.

One set  of applications of our  results is to the  question of regular
elements  and regular  eigenvalues  for reflection  cosets. The  vector
$v\in  V$ is  ($G$-)regular  if  $v$ does  not  lie  on any  reflecting
hyperplane of  $G$. The element $\gamma\in  \gamma G$ is regular  if it
has  a regular  eigenvector; if  $\gamma  v=\zeta v$,  then $\zeta$  is
called a  regular eigenvalue  and its order  (when $\gamma$  has finite
order) is a  regular number. In this work we  give precise criteria for
an eigenvalue  to be regular  for a coset,  and apply these  to various
questions. When $\gamma\in G$, it  is trivial that the identity element
of  $G$ is  a regular  element, with  corresponding eigenvalue  $1$. In
general, it  is not even  obvious that  any regular elements  exist. We
show that, with obviously necessary qualifications, they do.
\footnote{We have recently discovered that  this result also appears in
the work \cite{Ma} of G.~Malle, whom we thank for a preprint. Our proof
involves less case by case checking.}

In  the case  of ``well-generated  groups''  we use  our criterion  for
regularity which is  couched in terms of the $\eps_\iota$  to produce a
twisted analogue of Coxeter elements of  a real reflection group, and a
twisted analogue for a reflection coset of the Coxeter number of a real
reflection group.

Another significant application of our results concerning parabolic subgroups
is to the module structure of the coinvariant
algebra for  the group $\genby{G,\Gamma}$,  where $\Gamma$ is  a finite
subgroup of the normaliser of $G$ in $\GL(V)$.
Our  result, Theorem \ref{coinvariant induced}, generalises  one of  
Stembridge (in  the untwisted  regular
case),  whose  proof  goes  back   to  Springer's  computation  of  the
eigenvalues   of  a   regular  element   in  any   representation.  One
interpretation of  Stembridge's result is  that it gives  an expression
for the  $G$-module structure of  the sum of certain  graded components
of  the  coinvariant  algebra.  Our  result  generalises  this  in  two
ways; first,  by considering  a larger  class of  sums by  removing the
restriction  of  regularity, and  second,  by  considering the  twisted
structure of the coinvariant algebra. The statement in 
(\ref{coinvariant induced}) expresses the sum of certain graded 
components of the coinvariant algebra of $G$ as a representation induced 
from the coinvariant algebra of a parabolic subgroup. 

In the final section, we explore the relationship between our twisted
invariants and the reflection quotients of reflection groups studied in
\cite{BBR}. These are quotients $G/L$ of $G$ which act as reflection
groups on the tangent space at $0$ of $V/L$, where $L$ is a normal 
subgroup of $G$.  Here we are able to relate the constants and other
invariants of $G$ with those of its reflection quotients in the above 
sense.

\section{Background and Notation\label{section notation}}

Let $K$ be an algebraically closed field of characteristic $0$, and $V$ be
a $K$-vector space  of dimension $r$. Let $G\subset \GL(V)$  be a finite
subgroup  generated  by (pseudo)reflections.  Denote  by $\CA$  (or
$\CA(G)$ when appropriate) the corresponding set  of reflecting hyperplanes,
and for each $H\in\CA$ choose a  linear form $L_H\in V^*$ with kernel $H$. Let
$S$  be  the algebra  of  polynomial  functions  on  $V$; it may  be
identified with the symmetric algebra of the dual vector space $V^*$.
The  subalgebra $S^G$ of  $G$-invariant functions
is  a polynomial  algebra. Let $\CN$  be the
normaliser  of $G$  in  $\GL(V)$; this is  a (not necessarily  connected)
reductive group. Denote by $I$ the ideal of $S$ generated by elements
of $S^G$ with no constant term, and let $\CH$ be the space of $G$-harmonic
polynomials, i.e., the polynomials which are annihilated by all $G$-invariant
polynomial differential operators on $S$ with no constant term.
Then $\CH$ is $\CN$-stable and $I \oplus \CH = S$. So, by Chevalley's theorem,
the  natural  map  $S^G\otimes  \CH\to  S$  is  an
isomorphism  of  $\CN$-modules.  The  algebra $S_G=S/I$  is  called  the
{algebra of coinvariants}.  Again by a result of Chevalley, it
is isomorphic as a  $G$-module to $K[G]$, so that
$\CH$  is also  a $G$-submodule  of $S$,  isomorphic to  the regular
representation of $G$, which is stable under $\CN$.


Let  $M$  be  any finite dimensional  $G$-module;
then  $(S\otimes  M^*)^G\simeq  S^G\otimes
(\CH\otimes M^*)^G$ is free as $S^G$-module, and since $\CH$ realises the
regular  representation of  $G$,  it  is of rank  $m=\dim M$.
Notice that $S\otimes M^*$ is graded: we declare $\deg F\otimes y=\deg F$,
for $F$ a homogeneous element of $S$ and $y\in M^*$. If
$u_1,\ldots,u_m$  is a homogeneous linear basis of $(\CH\otimes M^*)^G$, it
is clearly an $S^G$-basis  of $(S\otimes  M^*)^G$. The numbers
$m_i(M):=\deg u_i$ are the $M$-exponents of $G$.

The following observations concerning $G$ and $\CN$ are useful, and will
often be used without comment.

\begin{remark}\label{propsgamma}
Since  $G$  is   generated  by  reflections,  there is  a  unique
decomposition  $V=V^G  \oplus  V_1  \oplus \ldots  \oplus  V_k$,  where
the  $V_i$   are  irreducible,  pairwise   non-isomorphic,  non-trivial
$G$-submodules  of $V$,  and  correspondingly $G=G_1\times\dots  \times
G_k$, where  $G_i$ acts  as an irreducible  reflection group  on $V_i$,
and  acts trivially  on the  other  summands. Of  course distinct  pairs
$(G_i,V_i)$ may be isomorphic as reflection groups.

Clearly $C_{\GL(V)}(G)  \simeq \GL(V^G)  \times K^\times  \times \ldots
\times  K^\times$   ($k$  copies   of  $K^\times$), while
 $\CN= \GL(V^G) \times  N_{\GL(V_1 \oplus
\ldots \oplus V_k)}(G_1 \times \ldots \times G_k)$. Moreover  $\CN  /
C_{\GL(V)}(G)$ is evidently a  finite group, since it acts faithfully  as a group
of automorphisms  of $G$. Thus if $\gamma \in
\CN$, there exists $n \ge 1$ such that $\gamma^n$ centralises $G$.

Now any element  $\gamma\in \CN$ is  of  the form  $x\gamma'$,
where $x\in\GL(V^G)\leq C_{\GL(V)}(G)$  and 
$\gamma'\in  N_{\GL(V_1 \oplus  \ldots \oplus
V_k)}(G_1 \times \ldots \times G_k)$.
Since $\gamma'$ permutes the subspaces $V_i$ and $\Hom_G(V_i,V_j)$ has
dimension at most $1$, there exists $z'\in K^\times\times \ldots\times K^\times$ such that
$\gamma'z'$ has finite order, and taking $z=xz'$, we see that
for any $\gamma\in\CN$ there is an element
$z \in  C_{\GL(V)}(G)$ such  that $z\gamma$  is of
finite order. If  $G$ is essential, i.e. if $V^G=0$,  then $x=1$ above,
so $z$ is semisimple, and $\CN$  is a
finite extension of a torus which centralises $G$; 
in this case, every element of $\CN$ is semisimple.
In general, since every unipotent element of  $\CN$
centralises $G$, the  action of any element $\gamma\in\CN$ on  $G$ coincides  with the
action of  its semisimple part on $G$. 
\end{remark}

Remark \ref{propsgamma} shows that the action on $G$ of an
arbitrary element $\gamma$ of $\GL(V)$ is induced by a 
semisimple element of finite order, and so no generality is lost by 
making the

\begin{quotation}
\noindent{\bf Hypothesis.  } {\it Henceforth  we take $\gamma$ to  be a
fixed semisimple element of $\GL(V)$, which normalises $G$.}
\end{quotation}

The  coset  $G\gamma$  will  be   referred  to  as  a  {\it  reflection
coset}. Let $M$  be a $\genby{G,\gamma}$-module on  which $\gamma$ acts
semisimply. We  shall define  some important constants  associated with
the coset $G\gamma$.

Since  $\gamma$  acts semisimply  on  $(\CH\otimes  M^*)^G$, the  basis
elements  $u_i$ above  may be  taken to  be eigenvectors  for $\gamma$.
For  each  such   $M$, denote by $\CB(M,\gamma)$  a  fixed homogeneous
basis   of  $(\CH  \otimes M^*)^G$  which consists of eigenvectors of
$\gamma$.  Given $\iota  \in \CB(M,\gamma)$,  denote by  $\eps_\iota(M)$,
or  $\eps_\iota$  when  unambiguous, the  corresponding  eigenvalue  of
$\gamma$,  and  by $m_\iota(M)=m_\iota$  the  degree  of $\iota$.  Thus
$\gamma\iota=\eps_\iota\iota$  for $\iota  \in \CB(M,\gamma)$,  and for
any $g\in  G$, $g\gamma\iota=\eps_\iota\iota$,  whence the  multiset of
pairs $\{(\eps_\iota,m_\iota)~|~\iota \in \CB(M,\gamma)\}$ depends only
on (the isomorphism  class of) $M$ and the coset  $G\gamma$, and not on
the choice of $\gamma\in G\gamma$ or on the basis $\CB(M,\gamma)$.

\begin{definition}\label{parameters} For  any $\genby{G,\gamma}$-module 
$M$,  the multiset  $\{\eps_\iota(M)= 
\eps_\iota\mid \iota \in \CB(M,\gamma)\}$ will be referred to as the
multiset of  $M$-factors of  $G\gamma$. 
\end{definition}             

\begin{remark}\label{changing gamma}
Let     $\zeta      \in     K^\times$and      let     $M$      be     a
$\genby{G,\gamma,\zeta    \Id_V}$-module    on    which    $\zeta\Id_V$
acts    as    multiplication    by    a    scalar,    say    $\zeta_M$.
Then   we   may    take   $\CB(M,\zeta^{-1}\gamma)=\CB(M,\gamma)$   and
we    have     $\eps_\iota(\zeta^{-1}\gamma)=\zeta_M    \zeta^{m_\iota}
\eps_\iota(\gamma)$ for every $\iota \in \CB(M,\gamma)$.
\end{remark}

The  cases  $M=V$  and  $M=V^*$   will  figure  prominently  below.  In
these  cases,  for  simplicity, we  write  $\CB(\gamma)=\CB(V,\gamma)$,
$\CB^*(\gamma)=\CB(V^*,\gamma)$.   When    $\gamma=\Id_V$,   we   write
$\CB(M)=\CB(M,1)$.

\begin{example}\label{derivation}
Let  $d  :  S  \to  S  \otimes  V^*$  be  the  unique  derivation
of $S$-modules  such
that  $d  X  =  1  \otimes  X$  for  every  element  $X  \in  V^*$.  If
$(X_1,\dots,X_r)$  is a  basis  of  $V^*$, then  for  $P\in  S$, $dP  =
\sum_{i=1}^r  \frac{\partial P}{\partial  X_i} \otimes  X_i$. Evidently
$d$ commutes with the action  of $\GL(V)$. For any $\BN$-graded algebra
$A$,  denote  by  $A_+$  the  (augmentation)  ideal  of  elements  with
no  degree  zero  term.  Since  $\CN$ is  reductive,  there  exists  an
$\CN$-stable graded subspace $Y$ of  $S^G$ such that $S^G_+ = (S^G_+)^2
\oplus  Y$.  Let  $(P_1,\dots,P_r)$  be a  homogeneous  basis  of  $Y$.
Then  it  is  well-known  that  the  natural  map  $S(Y)  \to  S^G$  is
an  $\CN$-equivariant  isomorphism  of algebras (here $S(Y)$ denotes the
symmetric algebra on $Y$),  i.e.  that  $S^G\cong
K[P_1,\dots,P_r]$, and  that $(dP_1,\dots,dP_r)$  is an  $S^G$-basis of
$(S  \otimes V^*)^G$.  Denote by  $\bar d$  the composite  map $\bar{d}
:  Y\overset{d}{\to}(S\otimes  V^*)^G  \overset{\eta}{\to}(\CH  \otimes
V^*)^G$ where $\eta:  S\otimes V^*\to \CH\otimes V^*$  is the extension
to $S\otimes V^*$ of the natural map $S \to S/I \simeq \CH$. Then $\bar
d$  is  an  isomorphism  of  degree $-1$  of  $\CN$-modules.  For  each
$\iota  \in \CB(\gamma)$,  let  $P_\iota  = \bar{d}^{-1}(\iota)$.  Then
$\{P_\iota~|~\iota \in  \CB(\gamma)\}$ is  another basis of  $Y$, which
consists  of  $\gamma$-eigenfunctions.  The  $P_\iota$ form  a  set  of
homogeneous basic $G$-invariants in $S$,  and since $\bar d$ has degree
$-1$,
\begin{equation}\label{degre exposants}
\deg P_\iota = m_\iota + 1
\end{equation}
where  $m_\iota$  are  the  usual  exponents  of  $G$.  Further,  since
$\bar{d}$ respects the action of $\CN$, we have
\begin{equation}\label{action gamma invariants}
\gamma(P_\iota)=\eps_\iota P_\iota.
\end{equation}
\end{example}

\begin{definition}\label{definition degrees}
For $\iota \in \CB(\gamma)$, write $d_\iota=m_\iota+1$.
The multiset of degrees of $G$ is
$\{d_\iota=\deg P_\iota~|~\iota \in \CB(\gamma)\}$.

Correspondingly,     for    $\iota     \in    \CB^*(\gamma)$,     write
$d_\iota^*=m_\iota-1$.  In  this  case  the $m_\iota$  are  called  the
coexponents of $G$, and the multiset of codegrees of $G$ is
$\{d_\iota^*~|~\iota \in \CB^*(\gamma)\}$.
\end{definition}

We shall be making use of the following result of Gutkin. In discussing
it, we  take $\gamma=1$ above, i.e.  the theorem will be  stated in the
``untwisted'' context. For any $G$-module  $M$ of finite dimension $m$,
let $N(M)=\sum_{\iota \in \CB(M)} m_\iota$. Given $H\in \CA$, denote by
$G_H$ the cyclic  (reflection) subgroup of $G$  comprising the elements
which  fix $H$  pointwise, and  set $N_H(M)=N(\Res^G_{G_H}  M)$. 
If $\xi_H$ is the unique non-trivial component character of the 
representation of $G_H$ on $V$, then any
character $\xi$ of $G_H$ is uniquely expressible as $\xi=\xi_H^e$, where 
$0\leq e\leq e_H-1$. Accordingly, if we write
$\Res^G_{G_H}  M^*=\oplus_{i=1}^m \xi_H^{e_i}$ with $0\leq e_i\leq e_H-1$,
then clearly $N_H(M)=\sum_{i=1}^m e_i$.
Observe  that  for any  $G$-module  $M$,
$S\otimes \Lambda M^*$  is a bigraded associative  algebra, where $\deg
(F\otimes x_1\wedge\dots\wedge x_j)=(i,j)$ for  $F\in S$ homogeneous of
degree $i$ and $x_1,\dots,x_j\in M^*$.  The following theorem is due to
Gutkin (\cf. \cite[2.10]{OS}):

\begin{theorem}[Gutkin]\label{gutkin}
Let $y_1,\ldots,y_m$  be a basis  of $M^*$. Then the  product
$\prod_{\iota \in \CB(M)} \iota$  in $S\otimes  \Lambda  M^*$ lies
in  $(S\otimes\Lambda^m M^*)^G$ and satisfies
$$\prod_{\iota \in \CB(M)} \iota \dotequal \prod_{H\in\CA}
L_H^{N_H(M)}\otimes  y_1\wedge y_2\ldots\wedge  y_m,$$
where $\dotequal$ denotes equality up  to multiplication by some
$\lambda\in K^\times$. In particular by comparing degrees, we
have $N(M)=\sum_{H\in\CA} N_H(M)$.
\end{theorem}

The polynomial $\prod_{H\in\CA}L_H^{N_H(M)}$ will be denoted by $\Gut_M$.

\begin{example}\label{gutkin examples}
If  $H  \in  \CA$,  let  $e_H=|G_H|$.  Then  since  $N_H(V)=e_H-1$  and 
$N_H(V^*)=1$,  we   get  $$\Gut_V  =  \prod_{H   \in  \CA}  L_H^{e_H-1} 
\quad\text{ and } \quad \Gut_{V^*} = \prod_{H \in \CA} L_H.$$           
\end{example}

We shall also require the next result, which is due  to  Orlik and Solomon
(\cf. \cite[3.1]{OS}).

\begin{theorem}[Orlik and Solomon]\label{OS} Let $M$ be a $G$-module of
dimension $m$. If $N(\Lambda^m(M))=N(M)$, then
$$(S\otimes\Lambda(M^*))^G\simeq S^G\otimes\Lambda\bigl( (\CH \otimes
M^*)^G\bigr).$$

Equivalently, in the above notation, $$(\CH\otimes\Lambda M^*)^G
\simeq \Lambda((\CH\otimes M^*)^G).$$
\end{theorem}

The next  lemma (\cf. \cite[pp  79-82]{OS}) shows that  Theorem \ref{OS}
can be applied to a certain class of representations of $G$ which
include the Galois conjugates of $V$.

\begin{lemma}\label{pregalois}
Suppose $M$ is any $G$-module in which the reflections of $G$ act as 
reflections. Then $N(\Lambda^m M)=N(M)$, where $m=\dim M$.
\end{lemma}

Since the Galois conjugates of $V$ clearly satisfy the conditions of
(\ref{pregalois}), an immediate consequence is

\begin{corollary}\label{galois}
If $\sigma \in \Gal(K/\BQ)$, then $N(\Lambda^r(V^\sigma))=N(V^\sigma)$.
\end{corollary}

For the convenience of the reader, and also since our proofs may be slightly 
more straightforward  than those in the literature, we  provide  
proofs of Theorem  \ref{gutkin} and  Lemma \ref{pregalois}  in Appendix  2 
below.                                                                  

\subsection{Some bilinear forms}\label{bilinear}
We complete this section by defining  some bilinear forms which will be
used extensively below. If $\Gamma$ is a subgroup of $\CN$ and $M$ is a
$\genby{G,\Gamma}$-module, the $S$-bilinear form  $(S \otimes M) \times
(S \otimes M^*) \to S$, given by
$(f \otimes x,f' \otimes \varphi) \mapsto \varphi(x) ff'$ is
$\genby{G,\Gamma}$-equivariant. Therefore it induces
by restriction an $S^G$-bilinear form
$\langle\;,\;\rangle_M : (S \otimes M)^G \times (S \otimes M^*)^G \to S^G$
which is $\Gamma$-equivariant. Take an element $\gamma \in \Gamma$.
For $\iota \in \CB(M,\gamma)$ and $\jmath \in \CB(M^*,\gamma)$, we set
$$\CM_{\iota \jmath}^M = \langle \iota, \jmath \rangle_M.$$
Evidently the matrix
$\CM^M=(\CM_{\iota \jmath}^M)_{(\iota,\jmath) \in
\CB(M,\gamma)\times\CB(M^*,\gamma)}$ has entries in $S^G$, and we have
\begin{equation}\label{gammamij}
\gamma(\CM_{\iota \jmath}^M) = \eps_\iota \eps_\jmath \CM_{\iota\jmath}^M.
\end{equation}
Let $\Delta_M\in S^G$ denote the determinant of $\CM^M$.
\begin{lemma}
We have
\begin{equation}\label{delta m}
\Delta_M = \Delta_{M^*} \dotequal \Gut_M \Gut_{M^*}.
\end{equation}
\end{lemma}
\begin{proof}
Let $(v_1,\ldots,v_m)$ be a basis of $M$ and $(v_1^*,\ldots,v_m^*)$
the dual basis of $M^*$. For 
$\iota\in\CB(M)$ write $\iota=\sum_{k=1}^m q_{\iota k}^* \otimes v_k^*$
and for $\jmath\in\CB(M^*)$  write
$\jmath=\sum_{k=1}^m q_{\jmath k} \otimes v_k$
where $q_{\jmath k},\;q^*_{\iota k}\in\CH\subset S$.
Let $Q=(q_{\jmath k})_{\jmath\in\CB(M^*), 1 \le k \le m}$ and
$Q^* = (q_{\iota k}^*)_{\iota \in\CB(M), 1 \le k \le m}$. Then
$\Gut_M \dotequal \det Q^*$
and $\Gut_{M^*} \dotequal \det Q$. Therefore,
$\Gut_M \Gut_{M^*} \dotequal (\det Q) (\det Q^*)$.
Now $\langle \iota, \jmath \rangle_M =
\sum_{k = 1}^m q_{\jmath k} q_{\iota k}^*: = r_{\jmath\iota}$,
where $Q~\lexp{t}{Q^*}=(r_{\jmath\iota})_{\jmath\in\CB(M^*),\iota\in\CB(M)}$.
Therefore $\Delta_M\dotequal\det(Q~\lexp{t}{Q^*})$, as stated.
\end{proof}

\begin{example}\label{discriminant}
Write  $\CM=\CM^V$  and $\Delta=\Delta_V$.  Then  $\CM$  is called  the
{\it  discriminant  matrix} of  $G$  and  its determinant  $\Delta$  is
the  {\it  discriminant}  of  $G$.  From the  above,  we  see  $$\Delta
\dotequal \prod_{H \in \CA} L_H^{e_H}$$  (see (\ref{delta m}) and Example
\ref{gutkin examples}).
\end{example}

\section{A Twisted Polynomial Identity\label{section identity}}

We start with the following ``twisted''  version of  a result of
Orlik and  Solomon (\cf. \cite[3.3]{OS}).

\begin{theorem}\label{OS2}
Let $M$ be a $\genby{G,\gamma}$-module of dimension $m$
such that $N(\Lambda^m(M))=N(M)$. Then
$$ |G|\inv\sum_{g\in G}\frac{\det(1-yg\gamma\mid M^*)}
{\det(1-xg\gamma\mid V^*)}
=\frac{\prod\limits_{\iota \in \CB(M,\gamma)}(1-y\eps_\iota x^{m_\iota})}
{\prod\limits_{\iota \in \CB(\gamma)} (1-\eps_\iota x^{d_\iota})}.$$
\end{theorem}
\begin{proof}
We   have   seen   that   $S\otimes  \Lambda   M^*$   is   a   bigraded
$K$-vector   space.   Thus   we   may  define   the   bi-graded   trace
$\Tr(\alpha;x,y)\in  K[[x,y]]$  of  a bi-graded  endomorphism  $\alpha$
by $$\Tr_{(S\otimes\Lambda  M^*)^G}(\alpha;x,y)=\sum_{i,j\geq 0}^\infty
\Tr(\alpha, (S\otimes\Lambda M^*)^G_{i,j})x^iy^j.$$

We  now  compute  $\Tr_{(S\otimes\Lambda  M^*)^G}(\gamma;x,y)$  in  two
different ways  using (\ref{OS}). On  the left side we  use a variant of
Molien's formula,  while on  the right, we  use well-known  methods for
computing graded traces in tensor and exterior algebras (cf. e.g., \cite{LM}).
\end{proof}

For  $\zeta\in K^\times$  and $g\in\GL(V)$  denote by  $V(g,\zeta)$ the
$\zeta$-eigenspace  of $g$.  Clearly  $V(g,\zeta)$  coincides with  the
subspace $V^{\zeta\inv g}$ of points of $V$ fixed by $\zeta\inv g$.

For  any   finite  dimensional  $\genby{G,\gamma}$-module   $M$,  write
$U(M,\gamma)$ for the  set $\{\iota \in \CB(M,\gamma)~|~\eps_\iota=1\}$
and   $U_{\#}(M,\gamma)=\CB(M,\gamma)   \setminus  U(M,\gamma)$.   Then
$U(M,\gamma)$    is   a    homogeneous   basis    of   $(\CH    \otimes
M^*)^{\genby{G,\gamma}}$.   In  particular,   $|U(M,\gamma)|=\dim  (\CH
\otimes M^*)^{\genby{G,\gamma}}$.

Since  $G$ is  finite, as  a $G$-module,  $V^*$ is  a Galois  conjugate
of  $V$.  However  this  is  not  the case  for  $V^*$  regarded  as  a
$\GL(V)$-module. Since  we include elements  $\gamma$ in our  discussion,
the inverses of whose eigenvalues may not be Galois conjugate 
(setwise) to those of $\gamma$, we
need to  distinguish between  the Galois
conjugates of $V$ and those of $V^*$.
For $\sigma \in \Gal(K,\BQ)$, write $\CB(\sigma,\gamma)$,
$U(\sigma,\gamma)$ and $U_\#(\sigma,\gamma)$ for
$\CB(V^\sigma,\gamma)$, $U(V^\sigma,\gamma)$ and
$U_\#(V^\sigma,\gamma)$ respectively.
Similarly, write $\CB^*(\sigma,\gamma)$,
$U^*(\sigma,\gamma)$ and $U_\#^*(\sigma,\gamma)$ for
$\CB((V^*)^\sigma,\gamma)$, $U((V^*)^\sigma,\gamma)$ and
$U_\#((V^*)^\sigma,\gamma)$ respectively.
Finally, write $U(\gamma)$,
$U_\#(\gamma)$, $U^*(\gamma)$ and $U_\#^*(\gamma)$ for
$U(V,\gamma)$, $U_\#(V,\gamma)$, $U(V^*,\gamma)$ and $U_\#(V^*,\gamma)$.
Thus, for example, $U^*(\gamma)=U(V^*,\gamma)$, which is a basis of
$(\CH\otimes V)^{\genby{G,\gamma}}$, and
$U(\gamma)=U(V,\gamma)$, which is a basis of
$(\CH\otimes V^*)^{\genby{G,\gamma}}$.

\begin{proposition}\label{inegalite}
For any $\gamma \in \CN$ and any
$\sigma \in \Gal(K/\BQ)$, we have~:
\begin{enumerate}
\item If $V \ne V^G$, then $|U^*(\gamma)| \ge 1$.

\item $|U(\gamma)| \le |U(\sigma,\gamma)|$ and
$|U(\gamma)| \le |U^*(\sigma,\gamma)|$.
\end{enumerate}
\end{proposition}

\begin{proof}
(i) If $(v_1\dots,v_r)$ and $(X_1,\dots,X_r)$ are dual bases of $V$ and 
$V^*$ respectively, the element $\sum_i  X_i\otimes v_i \in S\otimes V$ 
is invariant  under the whole of  $\GL(V)$, and hence a  fortiori under 
$\genby{G,\gamma}$. Moreover if $V^G=0$, this element lies in $\CH\otimes V$,
and so $\dim (\CH\otimes V)^{\genby{G,\gamma}} \geq 1$. More generally,
whenever $V\neq V^G$, it 
represents a non-zero invariant element of degree $1$ of $S_G\otimes V$, whence the 
statement. 

\smallskip

(ii) follows from the same argument as in \cite[Proof of Theorem 2.3]{LM},
applied to Theorem \ref{OS2} taking $M=V^\sigma$ and $M=(V^*)^\sigma$ respectively.
Note that both these choices of $M$ satisfy the condition of (\ref{OS2}) by
Lemma \ref{pregalois}.
\end{proof}

The  following  result   is  deduced  from  Theorem   \ref{OS2}  as  in
\cite[Theorem 2.3]{LM}. Note that (\ref{sigma}) is obtained by applying
Theorem  \ref{OS2} with $M=V^\sigma$  while (\ref{sigma  dual}) is  obtained
by  applying  Theorem  \ref{OS2} with  $M=(V^*)^\sigma$. Theorem \ref{OS2}
applies to both cases by Lemma \ref{pregalois}.

\begin{theorem}\label{gen}
If $h:  V\to V$ is  a linear  transformation, denote by  $\det'(h)$ the
product of the non-zero eigenvalues of  $h$. Then we have the following
polynomial identities in $K[T]$ for any $\sigma \in \Gal(K/\BQ)$. In the formulae
below $\det$ always refers to the determinant on $V$.
\begin{multline}\label{sigma}
\sum_{g\in G}T^{\dim V^{g\gamma}}{\textstyle\det'}
(1-g\gamma)^{\sigma-1}=\\
\begin{cases}
\quad 0 & \text{if $|U(\gamma)|\ne|U(\sigma,\gamma)|$,}\\
\displaystyle
\prod_{\iota \in U(\sigma,\gamma)}\!\!\!\!(T+m_\iota)
\prod_{\iota\in U_\#(\sigma,\gamma)} \!\!\!\!(1-\eps_\iota^{-1})
\prod_{\iota\in U_\#(\gamma)}\frac{d_\iota}{1-\eps_\iota^{-1}}\;
& \text{otherwise.}\\
\end{cases}
\end{multline}
\begin{multline}\label{sigma dual}
(-1)^r \sum_{g\in G}(-T)^{\dim V^{g\gamma}}{\textstyle\det'}
(1-g\gamma)^{\sigma-1}\det(g\gamma)^{-\sigma}=\\
\begin{cases}
\quad 0 & \text{if $|U(\gamma)|\ne|U^*(\sigma,\gamma)|$,}\\
\displaystyle
\prod_{\iota \in U^*(\sigma,\gamma)}\!\!\!\!(T+m_\iota)
\prod_{\iota\in U^*_\#(\sigma,\gamma)} \!\!\!\!\!\!(1-\eps_\iota^{-1})
\prod_{\iota\in U_\#(\gamma)}\!\frac{d_\iota}{1-\eps_\iota^{-1}}\;
& \text{otherwise.}\\
\end{cases}
\end{multline}
\end{theorem}

We record two special cases of this theorem.
They are obtained by taking $\sigma=\Id_K$ in (\ref{sigma}) and
(\ref{sigma dual}) respectively. Note that (\ref{LM2form}) shall be
reinterpreted in (\ref{better LM2form}) below.

\begin{equation}\label{twistpw}
\sum_{g\in G}T^{\dim V^{g\gamma}}=
\prod_{\iota\in U(\gamma)}(T+d_\iota-1)\prod_{\iota\in U_\#(\gamma)}d_\iota.
\end{equation}

\begin{multline}\label{LM2form}
(-1)^r\sum_{g\in G}\det(g\gamma)\inv(-T)^{\dim V^{g\gamma}}=\\
\begin{cases}
\quad 0 &\text{if }|U(\gamma)|\ne|U^*(\gamma)|,\\
{\displaystyle\prod_{\iota\in U^*(\gamma)}\!\!\!(T+d_\iota^*+1)
\prod_{\iota\in U_\#^*(\gamma)}\!\!\!(1-\eps_\iota^{-1})
\prod_{\iota\in U_\#(\gamma)}\!\!\frac{d_\iota}{1-\eps_\iota^{-1}}}
& \text{otherwise.}\\
\end{cases}\\
\end{multline}

We refer  to the  elements of $V-\bigcup_{H\in\CA}H$  as ($G$-)regular,
and call $\zeta\in \C^\times$ regular  for the coset $G\gamma$ if there
is  an  element of  $G\gamma$  which  has  a regular  eigenvector  with
corresponding eigenvalue  $\zeta$. In complete analogy  with \cite{LM},
we deduce the next statement from (\ref{LM2form}).

\begin{proposition}\label{oneregular} The eigenvalue $1\in K^\times$ is 
regular for  $G\gamma$ if  and only if  $|U(\gamma)|=|U^*(\gamma)|$, or 
equivalently                                                            
$\dim(\CH\otimes V^*)^{\genby{G,\gamma}}=\dim(\CH\otimes V)^{\genby{G,\gamma}}$.
\end{proposition}

\begin{remark}\label{changing gamma bis}
The element $\gamma\in\CN$ may be replaced by $\zeta^{-1}\gamma$, where
$\zeta$ is any element of $K^\times$, and the
formulae of Theorem \ref{gen} are then correspondingly modified, in a way
we shall now describe. Recall that as pointed out above,
$V^{\zeta^{-1} g \gamma}=V(g\gamma,\zeta)$. Further,
it follows from Remark \ref{changing gamma} that for any element
$\sigma \in \Gal(K/\BQ)$,
$\eps_\iota(\zeta^{-1}\gamma)=\zeta^{m_\iota + \sigma} \eps_\iota(\gamma)$
for each basis element $\iota \in \CB(\sigma,\gamma)$. Similarly,
$\eps_\iota(\zeta^{-1} \gamma) = \zeta^{m_\iota -\sigma} \eps_\iota(\gamma)$
for every $\iota \in \CB^*(\sigma,\gamma)$. Therefore,
from Definition \ref{definition degrees}, we have
$\eps_\iota(\zeta^{-1}\gamma) = \eps_\iota(\gamma)\zeta^{d_\iota}$ for every
$\iota \in \CB(\gamma)$ and we have
$\eps_\iota(\zeta^{-1}\gamma) = \eps_\iota(\gamma)\zeta^{d_\iota^*}$ for every
$\iota \in \CB^*(\gamma)$.
\end{remark}

Remark  \ref{changing  gamma  bis}  immediately  yields  the  following
general form of the  criterion (\ref{oneregular}) for regularity, which
is the twisted generalisation of the one given in \cite{LS2,LM}.

\begin{corollary}\label{oldregular}
The element $\zeta\in K^\times$ is regular for $G\gamma$ if and only if
$|\{\iota\in\CB(V,\gamma)\mid \eps_\iota\zeta^{d_\iota}=1\}|=
|\{\jmath\in\CB(V^*,\gamma)\mid \eps_\jmath\zeta^{d^*_\jmath}=1\}|.$
\end{corollary}

\section{Parabolic subgroups}

Let $v$ be any point in $V$ and let $C_G(v)=\{g \in G\mid g(v)= v\}$; this
is a  parabolic subgroup of $G$,  and contains as a normal subgroup the 
group $G_v$ defined as the subgroup of $C_G(v)$ which is generated by
reflections which fix $v$. Of course by Steinberg's Theorem (cf. e.g.
\cite{L}), the groups $G_v$ and $C_G(v)$ coincide, but we shall not assume 
this for the moment, since as a special case of the results of this section
we recover the proof of Steinberg's theorem, given in {\sl op. cit.}

Now $G_v$ is a reflection group and $\CA(G_v)=\{H  \in \CA(G)~|~v  \in H\}$.
Let $\Nv=  \{g \in\CN\mid g(v) \in K v\}$.
Then $\Nv$ contains the reflection group $G_v$ as a normal subgroup.

Let $I_v$ be the ideal of  $S$ generated by the homogeneous elements of
positive  degree  of $S^{G_v}$  and  write  $\CH_v$  for the  space  of
$G_v$-harmonic  polynomials,  i.e.  polynomials which  are  annihilated
by  all  $G_v$-invariant  polynomial  differential  operators  with  no
constant term. Evidently $\CH_v\subseteq \CH$, and $\Nv$ stabilises the
decomposition  $S=I_v  \oplus  \CH_v$; further, the  natural  map  $S^{G_v}
\otimes \CH_v \to S$ is an isomorphism of $\Nv$-modules. Notice that if
$v$ is $G$-regular, $G_v=\{1\}$, and $\CH_v=K$.

\begin{notation}
Let $N$ be any group and $M=\oplus_{i\in \bZ} M_i$ a $\bZ$-graded
$K[N]$-module. For any linear character
$\theta: N\to K^\times$ define $M\ogr \theta$
to be the graded $N$-module $\oplus_i M_i\otimes\theta^i$.
\end{notation}

Consider the linear map $\eta_v:S\to\CH_v$ given by
$f \otimes h  \in S\simeq S^{G_v} \otimes \CH_v \mapsto
f(v) h$. We also denote by $\eta_v  : \CH \to \CH_v$  its restriction
to $\CH$. Then $\eta_v$ clearly respects the action of $G_v$ on both sides;
we investigate how the action of $\Nv$ is transformed.
Let $\theta_v :  \Nv \to  K^\times$ be the  linear character
defined by $g(v) = \theta_v(g) v$ for every $g \in \Nv$.

\begin{lemma}\label{commnv} With the above notation,
$\eta_v$ induces an  epimorphism  of $\Nv$-modules from $\CH\ogr{\theta_v}$
to $\CH_v\ogr{\theta_v}$.
\end{lemma}
\begin{proof} Clearly $\eta_v$ is linear, and since $\CH\supseteq
\CH_v$, it is also evident that $\eta_v:\CH\to\CH_v$ is an epimorphism.
It therefore remains only to show that $\eta_v$ respects the indicated
actions of $\Nv$ on the two spaces.

For  $n\in\Nv$,  we  denote  simply  by $n$  its  action  on  $S$,  and
by  $\rho(n)$ (resp.  $\rho_v(n)$)  its  action on  $\CH\ogr{\theta_v}$
(resp. $\CH_v\ogr{\theta_v}$).  Then for  any element  $F=f\otimes h\in
(S^{G_v}\otimes \CH_v)\cap \CH$, with $f$ and $h$ homogeneous, we have
\begin{equation*}\begin{split}
\eta_v(\rho(n)(f\otimes h))&=\eta_v(\theta_v(n)^{\deg f+\deg h}
n(f)\otimes n(h))\\
&=\theta_v(n)^{\deg f+\deg h} n(f)(v) n(h)\\
&= \theta_v(n)^{\deg f+\deg h} f(n\inv(v)) n(h)\\
&= \theta_v(n)^{\deg f+\deg h} f(\theta_v(n)\inv v) n(h)\\
&= \theta_v(n)^{\deg h} f(v) n(h)\\
&= \rho_v(n)(f(v)h)\\
&= \rho_v(n)(\eta_v(f\otimes h)).\\
\end{split}
\end{equation*}
\end{proof}

Now let  $\Gamma$  be  any  subgroup  of   $\Nv$  and  let  $M$  be  a  finite
dimensional $\genby{G,\Gamma}$-module.
Consider $\CH\otimes M$  as a  graded module  by having
regard only to the degree in $\CH$ of its elements. Then the
$\Gamma$-modules $(\CH\ogr{\theta_v})\otimes M=(\CH\otimes M)\ogr{\theta_v}$
are canonically isomorphic. Observe that by Lemma \ref{commnv}, the map
$\eta_v^M:=\eta_v\otimes\Id: (\CH\ogr{\theta_v})  \otimes
M\to (\CH_v  \ogr{\theta_v})\otimes M$ respects the action of $\Gamma$
on both sides. Moreover $\eta_v^M$ clearly maps the $\Gamma$-submodule
$(\CH  \otimes M)^G\ogr{\theta_v}$ into
$(\CH_v \otimes M)^{G_v}\ogr{\theta_v}$.

\begin{theorem}\label{ced}
Let $\Gamma$ be any subgroup of $\Nv$ and let $M$ be a
finite dimensional $\genby{G,\Gamma}$-module.
Then the map $\eta_v^M$ introduced above induces  an isomorphism of
$\Gamma$-modules  from $(\CH \otimes M)^G\ogr{\theta_v}$ to
$(\CH_v \otimes M)^{G_v}\ogr {\theta_v}$.
\end{theorem}

\begin{proof}
Since $\eta_v^M$ is  $\Gamma$-equivariant by (\ref{commnv}),
it suffices to check that it is an
isomorphism of vector spaces. Let  $m=\dim M$ and let $y_1$,\dots, $y_m$
be a basis  of $M$. Let $u_1,\dots,u_m$  (resp. $u^v_1,\ldots,u^v_m$) be
homogeneous  bases of  $(\CH  \otimes  M)^G$ (resp.  $(\CH_v  \otimes
M)^{G_v}$).  Since $(\CH  \otimes  M)^G \subset  (S^{G_v} \otimes  \CH_v
\otimes M)^{G_v} =S^{G_v} \otimes (\CH_v \otimes M)^{G_v}$, we may write
$u_i=\sum_j q_{ji}u^v_j$ for some $q_{ji}\in S^{G_v}$.

We   now   apply   Gutkin's   Theorem  (see   (\ref{gutkin}))   in   turn
to   $G$    and   $G_v$;   since    $u_1\ldots   u_m=\det(q_{ij})_{i,j}
u^v_1\ldots     u^v_m$     we     obtain     $$\det(q_{ij})_{i,j}     =
\prod_{H\in\CA(G)-\CA(G_v)}L_H^{N_H(M)}.$$   But  for   every  $H   \in
\CA(G)-\CA(G_v)$, $L_H(v) \ne0$. Hence $\det(q_{ij})_{i,j}(v)\ne 0$.

Finally,  recall that  $\eta_v^M(u_i)=\sum_{j=1}^m q_{ji}(v)  u^v_j$ by
definition, whence $\eta_v^M$ is invertible.
\end{proof}

We may apply this result using a similar argument to that given in
\cite{L}, to relate the $M$-factors of $G$ and its parabolic subgroups.
Let $\zeta \in K^\times$ and let  $v\in  V$   be  a
$\zeta$-eigenvector of $\gamma$, so that $\gamma(v)=\zeta v$.
Let $M$ be a $\genby{G,\gamma}$-module.
Let $G_v$ be the stabiliser  of $v$ in
$G$; since  $\gamma\in N_{\GL(V)}(G_v)$, $M$ is also a
$\genby{G_v,\gamma}$-module, and we may  consider the
basis $\CB^v(M,\gamma)$ of $(\CH_v \otimes M^*)^{G_v}$ consisting
of homogeneous $\gamma$-eigenvectors. We also
define $U^v(M,\gamma)$ and $U_\#^v(M,\gamma)$ as analogues
for the pair $(G_v,\gamma)$ of the sets defined earlier for $(G,\gamma)$.


\begin{corollary}\label{eqlists}
Let $\gamma \in \Nv$, let $\zeta=\theta_v(\gamma)$ and let $M$ be a
$\genby{G,\gamma}$-module. Then the multisets
$\{\eps_\iota \zeta^{m_\iota}~|~\iota \in \CB(M,\gamma)\}$ and
$\{\eps_\iota \zeta^{m_\iota}~|~\iota \in \CB^v(M,\gamma)\}$
are equal.
\end{corollary}

\begin{proof}
Since $\theta_v(\gamma)=\zeta$,
$\{\eps_\iota \zeta^{m_\iota}~|~\iota \in \CB(M,\gamma)\}$
is the multiset of eigenvalues of $\gamma$ on
$(\CH \otimes M^*)^G \ogr \theta_v$ and
$\{\eps_\iota \zeta^{m_\iota}~|~\iota \in \CB^v(M,\gamma)\}$
is the multiset of eigenvalues of $\gamma$ on
$(\CH_v \otimes M^*)^{G_v} \ogr \theta_v$.
Applying  Theorem \ref{ced}, with $M$ replaced by $M^*$ and
$\Gamma$ by $\genby\gamma$, we obtain that these two multisets are equal.
\end{proof}

\begin{corollary}[Proof of Steinberg's theorem, cf. \cite{L}]
In the notation of the first paragraph of this section, we have
$C_G(v)=G_v$.
\end{corollary}
\begin{proof}
Take $\gamma\in C_G(v)$. It remains to show $\gamma\in G_v$. Since 
$\gamma v=v$ we take $\zeta=1$ and $M=V$ in (\ref{eqlists}), to obtain,
taking into account (\ref{action gamma invariants}),
that the multiset $\{\e_\iota~\mid~\iota \in \CB^v(M,\gamma)\}$ consists
entirely of $1$'s, since evidently $\gamma$ acts trivially on $S^G$.
Hence $\gamma$ acts trivially on $S^{G_v}$, whence $\gamma\in G_v$.
\end{proof}

\subsection{The coinvariant algebra as $\genby{G,\gamma}$-module}

We  shall  apply the  above  considerations  to  prove a  result  which 
generalises that  of Stembridge  \cite[2.3]{Ste} in  two ways.  Given a 
regular number $d$ for $G$, the result in {\sl loc. cit.} expresses the 
sum of the graded components of $S_G$ (or $\CH$) of degree congruent to 
$k$ modulo $d$ as an induced representation. Here we prove an analogous 
result without the  restriction that $d$ be regular;  we further extend 
the statement  to the action  of $\genby{G,\Gamma}$, where  $\Gamma$ is 
any finite subgroup of $\Nv$.       

Fix   such   a   subgroup   $\Gamma$  of   $\Nv$.   Note   that   since 
$\theta_v$   is  trivial   on   $G_v  \cap   \Gamma$,   it  defines   a 
linear  character  of   $\genby{G_v,\Gamma}$  through  the  isomorphism 
$\Gamma/(G_v  \cap  \Gamma)\simeq  \genby{G_v,\Gamma}/G_v$.  This  will 
also  be   referred  to   as  $\theta_v$.  It   thus  makes   sense  to 
consider   the  $\genby{G_v,\Gamma}$-module   $\CH_v\ogr\theta_v$,  and 
more  generally  for  any $k\in\BZ$,  the  $\genby{G_v,\Gamma}$-modules 
$(\CH_v\ogr\theta_v)\otimes\theta_v^k$. However, $\theta_v^k$ defines a 
character  of  $\genby{G,\Gamma}$  only  when $G  \cap  \Gamma  \subset 
\Ker\theta_v^k$.

\begin{theorem}\label{coinvariant induced}
Let $\CH_i$  denote  the  homogeneous component of degree $i$ of
the space $\CH$ of $G$-harmonic polynomials. Then maintaining the above
notation, for  any  integer $k\in\BZ$
there is an  isomorphism   of $\genby{G,\Gamma}$-modules
\begin{equation}\label{coineq}
\mathop{\oplus}_{\{i \ge 0\mid G \cap \Gamma \subset \Ker \theta_v^{k+i}\}}
\CH_i\otimes\theta_v^{k+i}\cong\Ind_{\genby{G_v,\Gamma}}^{\genby{G,\Gamma}}
\bigl((\CH_v\ogr\theta_v)\otimes\theta_v^k\bigr).
\end{equation}
\end{theorem}

\begin{proof}
It suffices to show that both sides have the same inner product with
any irreducible $\genby{G,\Gamma}$-module $M$. Let $X$ denote the
$\genby{G,\Gamma}$-module
on the left hand-side of the above equation.
Then $\langle X,M \rangle_{\genby{G,\Gamma}}=\dim ((X \otimes M^*)^G)^\Gamma$.
Therefore,
\begin{eqnarray*}
\langle X,M \rangle_{\genby{G,\Gamma}}
&=& \displaystyle{\frac{1}{|\Gamma|} \sum_{\gamma \in \Gamma}
\Tr(\gamma, (X \otimes M^*)^G)} \\
&=& \displaystyle{\frac{1}{|\Gamma|} \sum_{\gamma \in \Gamma}
\sum_{\iota \in \CB(M,\gamma)_k} \eps_\iota(\gamma) \theta_v(\gamma)^{m_\iota + k}},
\end{eqnarray*}
where  $\CB(M,\gamma)_k =  \{\iota  \in  \CB(M,\gamma)~|~G \cap  \Gamma 
\subset \Ker \theta_v^{m_\iota +  k}\}$. Let $[\Gamma/(G \cap \Gamma)]$ 
be a  set of  representatives of $\Gamma/(G  \cap \Gamma)$.  If $\gamma 
\in  [\Gamma/(G \cap  \Gamma)]$  and  $g \in  G  \cap  \Gamma$, we  can 
take $\CB(M,\gamma)=\CB(M,\gamma  g)$. We then  have $\eps_\iota(\gamma 
g)=\eps_\iota(\gamma)$ for every $\iota \in \CB(M,\gamma)$. Therefore   
\begin{eqnarray*}
&& \displaystyle{\frac{1}{|\Gamma|} \sum_{\gamma \in \Gamma}
\sum_{\iota \in \CB(M,\gamma)} \eps_\iota(\gamma) 
\theta_v(\gamma)^{m_\iota + k}}\\
&=&
\displaystyle{\frac{1}{|\Gamma|} \sum_{\gamma \in [\Gamma/(G \cap \Gamma)]}
\sum_{\iota \in \CB(M,\gamma)} \Bigl(\eps_\iota(\gamma) 
\theta_v(\gamma)^{m_\iota + k}
\sum_{g \in G \cap \Gamma} \theta_v^{m_\iota + k}(g) \Bigr)}\\
&=& \displaystyle{\frac{1}{|\Gamma|} \sum_{\gamma \in \Gamma}
\sum_{\iota \in \CB(M,\gamma)_k} \eps_\iota(\gamma) 
\theta_v(\gamma)^{m_\iota + k}}.
\end{eqnarray*}
It follows that
$$\langle X,M \rangle_{\genby{G,\Gamma}} =\frac{1}{|\Gamma|} 
\sum_{\gamma \in \Gamma}
\Bigl(\theta_v(\gamma)^k \sum_{\iota \in \CB(M,\gamma)}
\eps_\iota(\gamma) \theta_v(\gamma)^{m_\iota}\Bigr).
\leqno{(a)}$$
Now, let $X'$ be the $\genby{G,\Gamma}$-module on the right side
of (\ref{coineq}). By Frobenius reciprocity, we have
$$\langle X',M \rangle_{\genby{G,\Gamma}} =
\dim \bigl((\CH_v\ogr\theta_v)\otimes\theta_v^k 
\otimes M^*\bigr)^{\genby{G,\Gamma}}.$$
Therefore,
$$\langle X',M \rangle_{\genby{G,\Gamma}}
=\frac{1}{|\Gamma|} \sum_{\gamma \in \Gamma} \Bigl(\theta_v(\gamma)^k
\sum_{\iota \in \CB^v(M,\gamma)} \eps_\iota(\gamma) 
\theta_v(\gamma)^{m_\iota} \Bigr).
\leqno{(b)}$$
The result now follows from (a) and (b), given
Corollary \ref{eqlists}.
\end{proof}

\begin{remark}
The special  case when $\Gamma=\genby{\gamma}$ with
$\gamma\in G$ and $v$ is regular is in \cite[2.3]{Ste}.
In this  case (\ref{eqlists}) essentially amounts to
Springer's description of the eigenvalues of a regular element.
The general version above could in principle be used to determine
the individual graded components of $\CH$ as $\genby{G,\Gamma}$-modules.
\end{remark}

\begin{remark}
Maintain the notation of the previous theorem and
assume further that $\gamma \in G$; write $\zeta=\theta_v(\gamma)$,
and let $d$ be the order of $\zeta$. Then of course $\genby{G_v,\gamma}$
is contained in $G$
and Theorem \ref{coinvariant induced} can be written as follows.
$$\mathop{\oplus}_{i \equiv -k \text{ mod } d} \CH_i \cong
\Ind_{\genby{G_v,\gamma}}^{G}
\bigl((\CH_v\ogr\theta_v)\otimes\theta_v^k\bigr).$$
In the special case where $G\simeq {\mathfrak{S}}_n$ is the
symmetric group of degree $n$, $d \le n$, $\gamma$ is the
product of $[n/d]$ disjoint cycles of length $d$ and
$G_v \simeq {\mathfrak{S}}_{n-[n/d] d}$, we retrieve
a result of Morita and Nakajima \cite{morita}.

As a consequence of the previous equation, one obtains (for any $G$) that
\begin{equation}\label{eqdims}
\dim \bigl(\mathop{\oplus}_{i \equiv -k {\mathrm{~mod~}} d} \CH_i\bigr)=
\dim \bigl(\mathop{\oplus}_{i \equiv -l {\mathrm{~mod~}} d} \CH_i\bigr)
\end{equation}
for every $k$ and $l$ in $\BZ$ and any natural number $d$ which is the
order of an eigenvalue of some element of $G$ (i.e. which divides
some degree $d_i$ of $G$). This implies that $\frac{1-t^{d_i}}{1-t}$
divides the Poincar\'e polynomial of $\CH$, which of course is well
known.
\end{remark}

\section{Regularity}\label{regularity}

In this section we shall refine and provide a different approach to
the regularity result (\ref{oldregular}) above. As usual,
we refer to the elements of $V_\reg=V-\bigcup_{H\in\CA} H$ as ($G$-)regular,
and call $\zeta\in K^\times$
regular for  the coset  $G\gamma$ if  there is  an element  of $G\gamma$
which  has a  regular  eigenvector  with corresponding eigenvalue  $\zeta$.
An element of $G\gamma$ which has a regular eigenvector is called
($G\gamma$-)regular. Note that $\zeta$ is regular for the coset $G\gamma$ if
and only if $1$ is regular for the coset $G (\zeta^{-1} \gamma)$,
so that in the context of regularity for cosets, it suffices to consider
$1$-regularity.
In this section we shall give several criteria for
the coset $G\gamma$ to contain a regular element.

We start with properties of the quotient variety $V/G$ and
the action of $\gamma$ on it. The ring of regular functions on $V/G$ is
$K[V/G]=S^G$. If $J$ is a subset of $S^G$, we denote by $\CV(J)$
the closed subvariety of $V/G$ it defines. Let $I(\gamma)$
be the ideal of $S^G$ generated by $(P_\iota)_{\iota \in U_\#(\gamma)}$
(recall that the $P_\iota,\iota\in\CB(\gamma)$ form a set of basic
homogeneous invariants for $G$, (\cf. (\ref{derivation})),
$\gamma P_\iota=\eps_\iota P_\iota$, and $\iota \in U_\#(\gamma)
\iff \eps_\iota\ne 1$).
\begin{lemma}\label{fixed points}
We have $\CV(I(\gamma))=(V/G)^\gamma$.
\end{lemma}

\begin{proof}
Let $K^{\CB(\gamma)}\simeq\BA^r$ be the $K$-vector space of sequences
$(x_\iota)_{\iota \in \CB(\gamma)}$ of elements of $K$ indexed by $\CB(\gamma)$.
Then the map $\pi : V \to K^{\CB(\gamma)}$,
$v \mapsto (P_\iota(v))_{\iota \in \CB(\gamma)}$ is a morphism
of varieties (corresponding to the inclusion $S^G \hookrightarrow S$)
which induces an isomorphism $V/G \simeq K^{\CB(\gamma)}$.
If we endow $K^{\CB(\gamma)}$ with the linear action of $\gamma$ given by
$$\gamma.(x_\iota)_{\iota \in \CB(\gamma)} =
(\eps_\iota x_\iota)_{\iota \in \CB(\gamma)},$$
then, by (\ref{action gamma invariants}), the morphism $\pi$ is
$\gamma$-equivariant. But using obvious notation, the space
$(K^{\CB(\gamma)})^\gamma$
is naturally identified with $K^{U(\gamma)}$, and
the lemma follows.
\end{proof}

The variety $V_\reg/G$ has a convenient description in these terms.
Recall from Example \ref{discriminant}, that if $\Delta$ is the discriminant
polynomial of $G$, we have
\begin{equation}\label{equation regular}
V_\reg / G = V/G - \CV(\Delta).
\end{equation}

We  next   point  out  twisted   generalisations  of  the   results  of 
\cite{bessis}.   The   following  result   has   the   same  proof   as 
\cite[1.4,1.6]{bessis};   we  include   it   here   for  the   reader's 
convenience.                                                            

\begin{proposition}\label{bessis}
Let  $G\gamma$  be   a  reflection  coset and let $\zeta \in K^\times$.
Then
\begin{enumerate}
\item $\zeta$  is regular for  $G\gamma$ if and  only if
$\Delta \notin I(\zeta^{-1}\gamma)$.

\item If  $\zeta$ is  regular  for $G\gamma$ and
$U(\zeta^{-1}\gamma)=\{\iota_0\}$, then  $\Delta$ is  monic in $P_{\iota_0}$.
\end{enumerate}
\end{proposition}

\begin{proof}
(i) Clearly $\zeta$ is regular for $G\gamma$ if and only if $1$
is regular for $G\zeta\inv\gamma$.
Hence we may assume without loss, that $\zeta=1$.
But $1$ is regular for $G\gamma$ if and only if $(V_\reg/G)^\gamma$
is non-empty. By Lemma \ref{fixed points} and
(\ref{equation regular}), this is the case if and only
if $\CV(I(\gamma))$ is not contained in $\CV(\Delta)$, that is,
if and only if $\Delta$ is not in the radical of $I(\gamma)$.
The result follows because $I(\gamma)$ is clearly a radical ideal of $S^G$.

(ii) Given that $\zeta$ is regular for $G\gamma$ and that
$U(\zeta^{-1}\gamma)=\{\iota_0\}$,
it follows from (i) that $\Delta$ is non-zero modulo $I(\zeta^{-1}\gamma)$,
which is generated by $\{P_\iota\mid \iota\neq \iota_0\}$. So modulo $I(\zeta^{-1}\gamma)$,
$\Delta \equiv \lambda P_{\iota_0}^k$ for some $\lambda \in K^\times$ and
some $k \ge 1$. Since $\Delta$ is homogeneous, $\Delta$ is
monic in $P_{\iota_0}$.
\end{proof}

Proposition \ref{regular} below generalises \cite[Theorem 3.1 (ii)]{LM}
and is a more precise version of (\ref{oldregular}) above.
We shall require some preliminaries before proving it.
Let $\CN_v=\{n \in \CN~|~n(v)=v\}$. This is a normal subgroup
of $\Nv$ and $\Nv/\CN_v \simeq K^\times$. Let $\Gamma$ be a subgroup
of $\CN_v$. Let $M$ be a $\genby{G,\Gamma}$-module. Consider
the bilinear form
$$\begin{array}{cccc}
\langle \;,\; \rangle_M^v : & (\CH \otimes M)^G \times (\CH \otimes M^*)^G
& \longrightarrow & K \\
& (f,g) &\longmapsto & (\langle f,g\rangle_M)(v).
\end{array}$$
This is simply the evaluation at $v$
of the element $\langle f,g\rangle_M \in S^G$ (see Section
\ref{section notation}). Clearly, $\langle\;,\;\rangle_M^v$ is
$\Gamma$-invariant.

\begin{lemma}\label{dualite}
If $v$ is regular and $\Gamma \subset \CN_v$, then
$\langle \;,\; \rangle_v^M$ is a $\Gamma$-invariant perfect pairing.
\end{lemma}

\begin{proof}
Observe that the discriminant of the bilinear form $\langle \;,\; \rangle_M^v$
is equal to $\Delta_M(v)$. But by (\ref{delta m}), $\Delta_M(v)\ne0$
if $v$ is regular.
\end{proof}

\begin{proposition}\label{regular}
Let $\gamma \in \CN$. Then the following are equivalent:
\begin{enumerate}
\item[(1)] $1$ is regular for $G\gamma$.

\item[(2)] The multisets $\{\eps_\iota^{-1}~|~\iota \in \CB(\gamma)\}$
and $\{\eps_\iota~|~\iota \in \CB^*(\gamma)\}$ are equal.

\item[(3)] $|U(\gamma)|=|U^*(\gamma)|$.
\end{enumerate}
\end{proposition}

\begin{remark}
The equivalence of (1) and (3) follows from an argument
similar to \cite[Theorem 3.1 (ii)]{LM} (\cf. (\ref{oldregular}) above).
However the proof we provide here will not make use of the polynomial
identities stated in Section \ref{section identity}.
\end{remark}

\begin{proof}
(1) $\Rightarrow$ (2) Assume that $1$ is regular for $G\gamma$.
Let $v \in V$ be $G$-regular and such that $g\gamma(v) = v$
for some $g \in G$. Then, by Corollary \ref{dualite},
$\bigl((\CH \otimes V)^G\bigr)^*$ and $(\CH \otimes V^*)^G$
are isomorphic $\genby{g\gamma}$-modules, via the
perfect pairing $\langle \;,\; \rangle_v$.
Therefore, they are isomorphic as $\genby{\gamma}$-modules.
But $\{\eps_\iota^{-1}~|~\iota \in \CB(\gamma)\}$ is the multiset
of eigenvalues of $\gamma$ on $\bigl((\CH \otimes V)^G\bigr)^*$
and $\{\eps_\iota~|~\iota \in \CB^*(\gamma)\}$ is the multiset
of eigenvalues of $\gamma$ on $(\CH \otimes V^*)^G$.
The statement follows.

\smallskip

(2) $\Rightarrow$ (3) is trivial.

\smallskip

(3) $\Rightarrow$ (1) Assume that $|U(\gamma)|=|U^*(\gamma)|$.
By replacing $\gamma$ by $g_0 \gamma$ for some $g_0 \in G$, we may
assume that $\dim V^{g\gamma} \le \dim V^\gamma$
for every $g \in G$. We choose $v \in V^\gamma$ in ``general position'',
i.e. such that $G_v$ acts trivially on $V^\gamma$.
By Corollary \ref{eqlists}, the multisets
$\{\eps_\iota~|~\iota \in \CB(\gamma)\}$ and
$\{\eps_\iota~|~\iota \in \CB^v(\gamma)\}$ are equal, so
$|U(\gamma)|=|U^v(\gamma)|$. Similarly, $|U^*(\gamma)|=|U^{v*}(\gamma)|$.
This shows that $|U^v(\gamma)|=|U^{v*}(\gamma)|$.
We now have:
\begin{quotation}
\begin{quotation}
\noindent (a) $V^\gamma \subset V^{G_v}$;

\noindent (b) $\dim V^{g\gamma} \le \dim V^\gamma$ for every $g \in G_v$;

\noindent (c) $|U^v(\gamma)|=|U^{v*}(\gamma)|$.
\end{quotation}
\end{quotation}
We shall show that this implies that $G_v=1$. Note that (b) implies
that $\dim V^\gamma=|U^v(\gamma)|$. Let $V'$ be the unique
$G_v$-stable subspace of $V$ such that $V=V^{G_v} \oplus V'$. It is
$\gamma$-stable, and
the homogeneous component of degree $1$ of $\CH_v$ is
$V^{\prime *}$.
By (a), $(K \otimes V^\gamma) \oplus
(V^{\prime *} \otimes V)^{\genby{G_v,\gamma}}$ is contained in
$(\CH \otimes V)^{\genby{G_v,\gamma}}$. Therefore,
\begin{eqnarray*}
|U^{v*}(\gamma)|&\ge&\dim V^\gamma+
\dim {\mathrm{Hom}}_{\genby{G_v,\gamma}}(V',V) \\
&=&|U^v(\gamma)|+\dim {\mathrm{Hom}}_{\genby{G_v,\gamma}}(V',V).
\end{eqnarray*}
It follows from (c) that $\dim {\mathrm{Hom}}_{\genby{G_v,\gamma}}(V',V) = 0$.
Thus $V'=0$, and $G_v=1$ as required.
\end{proof}

\begin{remark}
The right-hand side  of (\ref{LM2form}) vanishes unless  $1$ is regular 
for $G\gamma$.  Hence it  may be simplified  as follows.  All $g\gamma$ 
which  contribute to  the highest  power  of $T$  on the  left side  of 
(\ref{LM2form})  are  conjugate,  and  by  \cite[6.4(v)]{Sp}  have  the 
same  determinant, which  is  equal to  $\prod_{\iota \in  \CB(\gamma)} 
\eps_\iota^{-1}$. This in turns yields a formula for                    
$\frac{\prod_{\iota\in U_\#^*(\zeta)}(1-\eps_\iota^{-1})}{\prod_{\iota \in
U_\#(\zeta)}(1-\eps_\iota^{-1})}$, which may be substituted into
(\ref{LM2form}). The result is
\begin{multline}\label{better LM2form}
\sum_{g\in G}\det(g\gamma)T^{\dim V^{g\gamma}}=\\
\begin{cases}
\quad 0 & \text{if }|U(\gamma)|\ne|U^*(\gamma)|,\\
{\displaystyle \prod_{\iota \in \CB(\gamma)} \eps_\iota^{-1}
\prod_{\iota\in U_\#^*(\gamma)}\!(T-d_\iota^*-1)
\prod_{\iota\in U_\#(\gamma)}d_\iota} &\text{otherwise.}\\
\end{cases}
\end{multline}
\end{remark}

Our final observation in this section is that if $\Delta$ is monic in some
basic invariant, then there is a natural regular number.

\begin{corollary}\label{monicreg}
Suppose that the discriminant $\Delta$ is monic in $P_{\iota_0}$
for some $\iota_0 \in \CB(\gamma)$. Let $\zeta \in K^\times$ be such
that $\zeta^{d_{\iota_0}} = \eps_{\iota_0}^{-1}$. Then $\zeta$ is regular
for $G\gamma$. In particular, the multisets
$\{\eps_\iota \zeta^{d_\iota} ~|~\iota \in \CB(\gamma)\}$
and $\{(\eps_\iota \zeta^{d_\iota^*})^{-1}~|~
\iota \in \CB^*(\gamma)\}$ are equal.
\end{corollary}

\begin{proof}
Note that $\iota_0 \in U(\zeta^{-1} \gamma)$ by Remark \ref{changing gamma}.
Therefore, by assumption, $\Delta$ does not belong to the ideal
$I(\zeta^{-1}\gamma)$. So, by Proposition \ref{bessis} (i),
$\zeta$ is regular for $G\gamma$. Now,
the last assertion follows from Proposition \ref{regular}
and from Remark \ref{changing gamma}.
\end{proof}

\section{A twisted generalisation of Coxeter elements}

In  this section  we focus  attention on  ``well-generated'' reflection 
groups. These include  the finite Coxeter groups,  the Shephard groups, 
i.e. symmetry groups  of regular polytopes, and some  others. To define 
them, we have the

\begin{observation}[Orlik and Solomon]\label{wellgen} Let $G$ be
an irreducible reflection group in $V$.
Suppose the degrees and codegrees of $G$ are ordered so that
$d_1\leq d_2\leq \dots\leq d_r$ and $d_1^*\geq d_2^*\geq\dots\geq d_r^*$.
Then the following two statements are equivalent.
\begin{enumerate}
\item $G$ is generated by $r=\dim V$ reflections.
\item We have $d_i+d_i^*=d_r$ for $i=1,2,\dots,r$.
\end{enumerate}
\end{observation}

The only (currently) known proof of (\ref{wellgen}) is empirical.
A reflection group satisfying the equivalent conditions of (\ref{wellgen})
is called {\em well-generated}.

Henceforth, we shall consistently
write $\CB(\gamma)=(\iota_1,\dots,\iota_r)$ and
$\CB^*(\gamma)=(\iota_1^*,\dots,\iota_r^*)$ and we set $d_{\iota_i}=d_i$ and
$d_{\iota_i^*}^*=d_i^*$. Write $P_i=P_{\iota_i}$, so that $\deg P_i=d_i$.
We also assume that this numbering satisfies $d_1 \le \dots \le d_r$ and
$d_1^* \ge \dots \ge d_r^*$. We also set $\eps_{\iota_i}=\eps_i$ and
$\eps_{\iota_i^*}=\eps_i^*$.
The next result is part of Bessis' \cite[Theorem 2.2]{bessis2}. Again
for convenience, we sketch a proof.

\begin{proposition}\label{wellgenmatrix}
Suppose $G$ is any irreducible reflection group which satisfies
 $d_i+d_i^*\leq d_r$ for all $i$. Then
\begin{enumerate}
\item Any primitive $d_r^{\text{th}}$ root of unity is regular for $G$.
\item We have $0<d_i+d_j^*<2d_r$ for any pair $(i,j)$.
\item Let $I_0$ be the ideal of $S^G$ generated by $\{P_i\mid i\neq r\}$.
Then modulo $I_0$, the discriminant matrix $M\equiv P_rC$, where
$C=(c_{ij})_{i,j}$ is a non-singular matrix with entries in $K$.
\item If $d_i+d_j^*\neq d_r$, $c_{ij}=0$.
\item We have $rd_r=N+N^*$ and $d_i + d_i^*=d_r$ for every $i$.
\item Partition $\{1,\dots,r\}$ into subsets, where $i,j$ are in the same
subset if $d_i=d_j$. Then $C$ is diagonal by block for this decomposition
and each block is non-singular.
\item $G$ is well-generated.
\end{enumerate}
\end{proposition}
\begin{proof}
Since $G$ is irreducible, we have $d_i^* \ge 1$ for every $i \le r-1$.
Therefore, by assumption $1 \le d_i \le d_r-1$ and $1 \le d_i^* \le d_r-1$
if $i\ne r$.
So, if $\zeta_0$ is a primitive $d_r^{\text{th}}$ root of unity, then
$U(\zeta_0^{-1}\Id_V)=\{\iota_r\}$ and $U^*(\zeta_0^{-1}\Id_V)=\{\iota_r^*\}$,
whence $\zeta_0$ is regular for $G$ by
the criterion (\ref{regular}), proving (i); (ii) is a simple
consequence of our assumption on the degrees. By (\ref{bessis}),
$\Delta \equiv cP_r^k \mod I(\zeta_0^{-1} \Id_V)$, with $c\in K^\times$ by regularity,
and by degree, $k=(N(V)+N(V^*))/d_r$. Now the entries of $\CM$ are
homogeneous polynomials in $S^G$, and by (ii), when written as
polynomials in $P_r$ with coefficients in $K[\{P_i\mid i\neq r\}]$,
have degree (in $P_r$) $0$ or $1$. The statements (iii), (iv) and 
the first statement of (v) follow immediately. But,
$N(V)+N(V^*)=\sum_{i=1}^r (d_i+d_i^*)  \le rd_r = N(V)+N(V^*)$. So
$d_i+d_i^*=d_r$ for every $i$. This proves (v).
(vi) now follows from (iv) and (v) while (vii) follows from
(\ref{wellgen}).
\end{proof}

The next result refines (\ref{monicreg})  in the case of well-generated 
groups. It  may be regarded as  a generalisation of the  fact that when 
$G$  is real and crystallographic,  the  Coxeter number  $h$ (which  is
the highest  degree 
$d_r$), and  any primitive  $h^{\text{th}}$ root  of unity  is regular, 
with corresponding regular conjugacy class the Coxeter class of $G$.    

\begin{proposition}\label{well}
Suppose that $G$ is irreducible and well-generated. Let $\zeta \in K^\times$.
Then:
\begin{enumerate}
\item there is a permutation $\sigma$
of $\{1,\dots,r\}$ such that for each $i$, $d_i=d_{\sigma(i)}$
and $\eps_{\sigma(i)}^*=\eps_i^{-1}\eps_r$.
\item If $\zeta \in K^\times$ is such that $\zeta^{d_r}=\eps_r^{-1}$,
then $\zeta$ is regular for the coset $G\gamma$.
\end{enumerate}
\end{proposition}
\begin{proof}
By Proposition \ref{wellgenmatrix} (vi), there is a permutation $\sigma$ of
$\{1,\dots,r\}$ such that $d_{\sigma(i)}=d_i$ and
$c_{i,\sigma(i)}\neq 0$ for all $i$. Hence $\CM_{i,\sigma(i)}=
c_{i,\sigma(i)}P_r +\text{ other terms}$, from which it follows
that $\gamma(\CM_{i,\sigma(i)})=\eps_r \CM_{i,\sigma(i)}$. But
by (\ref{gammamij}), $\gamma(\CM_{i,\sigma(i)})=\eps_i\eps_{\sigma(i)}^*
\CM_{i,\sigma(i)}$.

Since $G$ is well-generated, we have $d_i+d_i^*=d_r$. Therefore,
by (i), $\eps_{\sigma(i)}^*\zeta^{d_{\sigma(i)}^*}=\eps_i^{-1}\zeta^{-d_i}$.
In view of the criterion (\ref{regular}), (ii) follows from this observation
(see also Remark \ref{changing gamma bis}).
\end{proof}

Finally, observe that  when $G$ is  real (and hence is a finite
Coxeter  group)  we  have, using our orderings,  $d_i^*=d_{r+1-i}-2$
and   $\eps_i^*=\eps_{r+1-i}$.  Thus, applying (i), we deduce that
there   is   a   degree-preserving  permutation   $\sigma$   such   that
$\eps_i\eps_{r+1-\sigma(i)} =\eps_r$.

\section{Existence of regular elements in cosets}\label{existence}

We shall prove
\footnote{As mentioned in the footnote to the Introduction, Theorem
\ref{regexists} also appears, with a different proof, in \cite{Ma}.}

\begin{theorem}\label{regexists}
There is a semisimple element
$z\in \GL(V)$ which centralises $\langle G,\gamma\rangle$
such that the reflection coset $z\gamma G$  has a regular
eigenvalue (or element).
\end{theorem}

As an easy consequence, we have

\begin{corollary}\label{reginter}
If $V$ is irreducible as $\genby{G,\gamma}$-module, then $G\gamma$
has a regular eigenvalue.
\end{corollary}

In particular,
\begin{corollary}\label{regirred}
If $V$ is irreducible as $G$-module, then $G\gamma$
has a regular eigenvalue.
\end{corollary}

We begin with a reduction to the case (\ref{regirred}),
which involves arguments similar to those in \cite[Prop. 6.9]{BL}.

\begin{lemma}\label{reduction}
We have the implications (\ref{regirred})$\implies$
(\ref{reginter})$\implies$(\ref{regexists}).
\end{lemma}

\begin{proof} To  see that (\ref{reginter})$\implies$(\ref{regexists}), 
suppose   that   $V=\oplus_i   V_i$   is   a   decomposition   of   $V$ 
into  irreducible  $\genby{G,\gamma}$-submodules. Then  correspondingly 
$G=G_1\times  G_2\times\dots$,  and $\gamma=\oplus_i  \gamma_i$,  where 
$G_i$  acts   as  a  reflection   group  in  $V_i$  and   trivially  on 
$V_j$  for  $j\neq i$,  and  $\gamma_i\in  \GL(V_i)$ normalises  $G_i$. 
The  set  $\CA$   of  reflecting  hyperplanes  of  $G$   is  the  union 
of  the  sets  $\CA_i$  of  reflecting hyperplanes  of  the  $G_i$.  By 
(\ref{reginter})  there  are elements  $g_i\in  G_i$  and $v_i\in  V_i$ 
such that  $\gamma_ig_iv_i=\zeta_iv_i$, and  $v_i$ is  $G_i$-regular in 
$V_i$. Take  $z=\oplus_i\zeta_i\inv\Id_{V_i}$, $g=(g_1,g_2,\dots)$, and 
$v=\oplus_iv_i$. Then $v$  is $G$-regular and $z\gamma  g v=v$, proving 
(\ref{regexists}).                                                      

Now    assume   (\ref{regirred}),    and    suppose    that   $V$    is 
irreducible   as   $\genby{G,\gamma}$-module.   Then   as   in   Remark 
\ref{propsgamma},  $V=V_1\oplus\dots  \oplus V_k$  and  correspondingly 
$G=G_1\times\dots\times  G_k$.  Then  all $(G_i,V_i)$  are  isomorphic, 
and   are   permuted   cyclically    by   $\gamma$.   Thus   $\gamma^k$ 
fixes   all   the   $V_i$,   and   in   particular   normalises   $G_1$ 
on   $V_1$,   so  that   by   (\ref{regirred}),   there  are   elements 
$g_1\in  G_1$   and  $v_1$  regular   in  $V_1$  such   that  $\gamma^k 
g_1v_1=\zeta_1v_1$.  Let  $\zeta\in \C$  satisfy  $\zeta^{-k}=\zeta_1$, 
let    $g=(1,1,\dots,1,\gamma^{k-1}g_1\gamma^{-(k-1)})\in    G$,    and 
$v=v_1\oplus\zeta\gamma  v_1\oplus  (\zeta\gamma)^2v_1\oplus\dots\oplus 
(\zeta\gamma)^{k-1}v_1\in  V$. Then  $v$  is  $G$-regular, and  $\gamma 
gv=\zeta\inv v$. Hence (\ref{regirred})$\implies$(\ref{reginter}).      
\end{proof}
It  follows  from Lemma  \ref{reduction},  that  it suffices  to  prove 
(\ref{regirred}),  and  hence   we  take  $V$  to   be  an  irreducible 
$G$-module. The next lemma deals with an obvious case.                  

\begin{lemma}\label{inner}
If $G$ is irreducible and $\gamma$ induces an inner automorphism
of $G$, then $($\ref{regirred}$)$ holds.
\end{lemma}

\begin{proof}
By hypothesis, there exists $g \in G$ such that $g \gamma$ is central
in $\genby{G,\gamma}$ so is scalar. The result is then obvious.
\end{proof}

At this point, it would be sufficient to inspect the list of
cosets $G\gamma$ such that $G$ is irreducible and $\gamma$ induces
a non-inner automorphism of $G$ which is given in \cite[3.13]{spets}.
The table of regular eigenvalues are then given in the table at the end
of this paper. However, we will provide some further reductions which
cover all the cases to be checked. First, we reduce the  proof further
to the case of  ``minimal groups'',which are defined as follows.

For any integer $d$, let $\zeta_d$ be a primitive $d^{\text{th}}$
root of unity. Then (\cf. \cite{LS1, LS2})
all maximal $\zeta_d$ eigenspaces $E$ of elements
of $G$ are conjugate under $G$, and the group $G(d):=N_G(E)/C_G(E)$
is a reflection group in $E$; the subquotient $G(d)$ is unique up to
conjugacy in $G$, and is irreducible if $G$ is \cite{LS2}. The regular case
is when $C_G(E)=1$, in which case $G(d)$ is a subgroup of $G$. Say that
$G$ is {\em minimal} if $\dim V>1$ and
there is no non-trivial subgroup $G(d)<G$ with $d$
regular. Equivalently, if $a(d),b(d)$ respectively denote the number
of degrees and codegrees divisible by $d$, then $a(d)=b(d)$ implies
that $a(d)=0$ or $r(=\dim V)$. Note that if $\dim V=1$, (\ref{regirred})
is trivially true.

\begin{lemma}\label{redmin}
Theorem \ref{regexists} is true for irreducible $G$
if it is true for irreducible minimal $G$.
\end{lemma}

\begin{proof}
If $\gamma$ normalises $G$ and  $E=V(g,\zeta_d)$ is a maximal $\zeta_d$
eigenspace,  then $\gamma  E=V(\gamma g\gamma\inv,\zeta_d)$  is also  a
maximal $\zeta_d$ eigenspace, whence there  is an element $x\in G$ such
that $\gamma E=x E$, so that $x\inv \gamma$ normalises $G(d)$, which is
irreducible by \cite[Theorem  A]{LS2}. If $v\in E$  is a $G(d)$-regular
eigenvector  for  $\gamma  y\in  x\inv \gamma  G(d)$,  then  since  the
reflecting  hyperplanes of  $G(d)$ are  the intersections  with $E$  of
those of  $G$, and $E$  is not contained in  any hyperplane of  $G$, it
follows that $v$ is $G$-regular. Thus Theorem \ref{regexists} holds for
$G$ if  it holds for  $G(d)$. Repeating this  argument, we arrive  at a
case where $G(d)$ is minimal.
\end{proof}

Our final lemma treats a case which arises frequently.

\begin{lemma}\label{ref}
Assume that $G$ is irreducible and that there exists $\gamma_0 \in \CN$ such
that $\genby{G,\gamma} \subset \genby{G,\gamma_0}$ and $\genby{G,\gamma_0}$
is a well-generated finite reflection subgroup of $\GL(V)$. Then
$G\gamma$ has a regular eigenvalue.
\end{lemma}

\begin{proof}
If $G$ is well-generated, the conclusion follows from Proposition
\ref{well} (ii). Hence we assume that $G$ is not well-generated.
By hypothesis, there exists $k \in \BZ$ such that $G\gamma_0^k=G\gamma$.
So if $\zeta$ is a regular eigenvalue for $G\gamma_0$, then $\zeta^k$
is a regular eigenvalue for $G\gamma$. Hence we are reduced to the case
$\gamma=\gamma_0$. Now let $Y$ be as in Example \ref{derivation}.
Then $S(V^*)^{\genby{G,\gamma}} \simeq S(Y)^\gamma$ is a polynomial algebra,
so $\gamma$ acts on $Y$ as a reflection. Let $i_0$ be the unique
element of $\{1,2,\dots,r\}$ such that $\eps_{i_0} \ne 1$ and let
$e$ be the order of $\eps_{i_0}$. Write $\widetilde{G}=\genby{G,\gamma}$.
Let $\tilde{d}_1 \le \dots \le \tilde{d}_r$ be the degrees of $\widetilde{G}$
and let $\tilde{d}_1^* \le \dots \le \tilde{d}_r^*$ be its codegrees.
Since $|U(\gamma)|=r-1$, two cases may occur
(see Proposition \ref{inegalite} (ii)).

If $|U^*(\gamma)|=r-1$, then $\gamma$ is $1$-regular.

If $|U^*(\gamma)|=r$, this means that $\gamma$ acts trivially on
$(\CH \otimes V)^G$. In particular, $\tilde{d}_i^*=d_i^*$ for every
$i$. Also, $S(Y)^\gamma\simeq S(V^*)^{\widetilde{G}}$
is a polynomial algebra generated $P_1$,\dots, $P_{i_0-1}$, $P_{i_0}^e$,
$P_{i_0+1}$,\dots, $P_r$. Therefore, since $\widetilde{G}$ is well-generated
and $G$ is not well-generated, it follows from Proposition \ref{wellgenmatrix}
that $(\tilde{d_1},\dots,\tilde{d}_r)=
(d_1,\dots,d_{i_0-1},d_{i_0+1},\dots,d_r,ed_{i_0})$.
Now, let $\zeta$ be such that $\zeta^{d_{i_0}}=\eps_{i_0}^{-1}$.
Then, since $\tilde{d}_i+\tilde{d}_i^* = ed_{i_0}$ and
$\zeta^{ed_{i_0}}=1$, $\eps_i=1$ if $i\ne i_0$, and that $\eps_i^*=1$
for every $i$ and since we have:

$\bullet$ If $1 \le i \le i_0-1$, then $d_i+d_i^*=ed_{i_0}$ and
$(\eps_i \zeta^{d_i})^{-1} = \eps_i^* \zeta^{d_i^*}$.

$\bullet$ $\eps_{i_0} \zeta^{d_{i_0}}=1=\eps_r \zeta^{d_r^*}$.

$\bullet$ If $i_0 + 1 \le i \le r$, then $d_i + d_{i-1}^* = ed_{i_0}$ and
$(\eps_i \zeta^{d_i})^{-1} = \eps_{i-1}^* \zeta^{d_{i-1}^*}$.

Therefore, the multisets
$\{\eps_i \zeta^{d_i}~|~1 \le i\le r\}$ and
$\{(\eps_i^* \zeta^{d_i^*})^{-1}~|~1 \le i\le r\}$ are equal.
So $\zeta$ is a regular eigenvalue for $G \gamma$
by Proposition \ref{regular}.
\end{proof}
We are now able to give the

\begin{proof}[Proof of Theorem \ref{regexists}]
The list of cosets $G\gamma$ such that $G$ is irreducible and $\gamma$
induces a non-inner automorphism of $G$ is given in \cite[3.13]{spets}.
Among them, the minimal ones are (up to multiplication by a scalar):

\begin{quotation}
\noindent{\bf 1.}
$G(de,e,r)\gamma$ when $r|e$, $d>1$ and $e>1$ where $\gamma \in G(de,1,r)$.

\smallskip

\noindent{\bf 2.} $G(4,2,2)\gamma$ where $\genby{G(4,2,2),\gamma}= G_6$.

\smallskip

\noindent{\bf 3.} $G_7 \gamma$ where $\genby{G_7,\gamma}=G_{15}$.
\end{quotation}

\noindent Since $G(de,1,r)$, $G_6$ and $G_{15}$ are well-generated,
these cases are disposed of by Lemma \ref{ref}.
The proof of the theorem is now completed by invoking
Lemmas \ref{reduction}, \ref{inner}, and \ref{redmin}.
\end{proof}


\section{Reflection quotients of reflection groups}

Let $L$ be a normal subgroup of $G$ and denote by $\bar{G}=G/L$ the
corresponding quotient. For any $\bZ_{\ge 0}$-graded algebra $A$,
denote by $A_+$ the ideal $\oplus_{i\geq 1}A_i$. Let 
$E^*$ be a graded complement of $(S_+^L)^2$ in $S_+^L$, so that 
$S_+^L=E^* \oplus (S_+^L)^2$. Evidently $E^*$ has basis a set of homogeneous 
generators of the invariant ring $S^L$. Let $\CN_L$ be the normalizer of $L$ 
in $\GL(V)$. For this section only, we denote the normalizer of
$G$ in $\GL(V)$ by $\CN_G$. Since $\CN_L$ 
is a reductive group, we may assume that $E^*$ 
is chosen to be stable under the action of $\CN_L$. 
Let $E$ be the (graded) dual of $E^*$, and denote 
by $\bar{S}$ the symmetric algebra of $E^*$. 

Then $E$ is isomorphic to the tangent 
space of the variety $V/L$ at $0$.
The quotient $\bar G$ acts on $E$, and we shall be interested in this
section in the case where this is a reflection group action, a situation which
has been studied in \cite{BBR}. In that case
$\CN_G\cap\CN_L$ also acts on $E$, normalising the $\bar G$-action,
and we shall relate the various twisted invariants of reflection cosets
of $G$ and $\bar G$.

The algebra homomorphism $\tau : S(E^*) \to S^L$ which extends the 
inclusion $E^*\hookrightarrow S^L$ is easily seen to be surjective
(see for instance \cite[Lemma 2.1]{BBR}) and $\CN_L$-equivariant. 
Denote by $I$ its ($\CN_L$-stable) kernel. Then 
we have a commutative diagram
$$\diagram
 & && && S & \\
0 \rto & I \rrto && \bar{S} \rrto^{\displaystyle{\tau}} && S^L \uto|<\ahook \rto & 0 \\
0 \rto & I^{\bar{G}} \rrto \uto|<\ahook && 
\bar{S}^{\bar{G}} \rrto^{\displaystyle{\tau^G}} \uto|<\ahook 
&& S^G \uto|<\ahook \rto & 0
\enddiagram$$
in which the rows are exact.

Note that the surjective 
morphisms $\tau$ and $\tau^G$ induce closed immersions

\begin{equation}\label{immersions}
V/L \hookrightarrow E\text{  and  } V/G \hookrightarrow E/\bar{G}.  
\end{equation}

We assume henceforth that $\bar{G}$ acts on $E$ as a reflection group. 
By \cite[Theorem 3.2]{BBR}, this is equivalent to the following 
two requirements:
\begin{quotation}
\noindent(1) $V/L$ is a complete intersection variety;

\noindent(2) $I$ is generated by $I^{\bar{G}}$ (as an ideal of $\bar{S}$).
\end{quotation}

\begin{example}
The classification of all such pairs $(G,L)$ is given in 
\cite[\S 4]{BBR}. If $L$ is generated by reflections, 
then it is straightforward that $G/L$ is always generated by reflections.
We give here another 
example (for more details, see \cite[\S 4.8]{BBR}). 
Assume that $G = G_{31}$ and that $L$ is its maximal normal 
$2$-subgroup. Then $L \subset {\text{SL}}(V)$, $|L|=64$, $|G/L|=720$, 
$E$ is of dimension $5$ and concentrated in degree $4$,
since $S^L$ is generated by $5$ polynomials of degree $4$. Here 
$G/L$ acts faithfully on $E$ as a group generated 
by reflections, isomorphic to the symmetric group ${\mathfrak{S}}_6$
in its irreducible reflection representation. Note that $G/L$ 
is well-generated while $G=G_{31}$ is not. 
\end{example}

\subsection{Coinvariants}
Denote by $\bar{S}_{\bar{G}}$ the algebra of coinvariants of $\bar{G}$. 
We shall relate this algebra to the algebra $S_G$ of coinvariants of $G$.

\begin{proposition}\label{tau iso}
The homomorphism  
$\tau:\bar S\lr S^L$ introduced above induces an isomorphism of graded algebras 
$$\bar{\tau} : \bar{S}_{\bar{G}} \longrightarrow (S_G)^L$$
which commutes with the action of $\CN_G \cap \CN_L(\supseteq G)$. 
\end{proposition}

\begin{proof} The composite of $\tau$ with the inclusion $S^L\hookrightarrow S$
maps $\bar S$ to $S$, and $\bar S_+^{\bar G}$ to $S_+^G$. Hence it induces
a homomorphism
$\bar \tau:\bar S_{\bar G}\to S_G$, whose image is evidently in $(S_G)^L$.
The equivariance with respect to $\CN_G\cap\CN_L$ is clear.

To prove that $\bar \tau$ is an isomorphism,
first note that $\bar{\tau}$ is surjective since $\tau$ is, because $L$ 
acts semisimply on $S$. 
But $\dim \bar{S}_{\bar{G}} = |\bar{G}|=\dim (S_G)^L$, whence $\bar\tau$
is also injective. 
\end{proof}

Let $\bar{\CH}$ be the space of $\bar{G}$-harmonic polynomial 
functions on $E$, and as above, $\CH$ be the corresponding space for
$G$ on $V$.

\begin{corollary}\label{harmonic} The isomorphism $\bar \tau$ of 
(\ref{tau iso}) induces an isomorphism of $\CN_G\cap\CN_L$-spaces
$:\bar{\CH} \lr \CH^L$, which we shall also denote by $\bar\tau$.
\end{corollary}

\begin{proof}
Each coset of the ideal $S. S^G_+$ of $S$ contains a unique $G$-harmonic
polynomial. This provides a canonical $\CN_G$-equivariant isomorphism
of vector spaces $:S_G\lr \CH$. Similarly we have a canonical 
$\CN_G\cap\CN_L$-equivariant canonical isomorphism $:\bar\CH\lr S_{\bar G}$.
If we compose $\bar\tau$ with these isomorphisms, taking (\ref{tau iso}) into 
account, we obtain the desired isomorphism.
\end{proof}

\subsection{Comparison of ${\boldsymbol{M}}$-factors}
Let $\Gamma$ be a subgroup of $\CN_G \cap \CN_L$, and write 
$\bar{\Gamma}\cong \Gamma/(\Gamma \cap L)$ for its image in
$\GL(E)$. 
Then $\bar{\Gamma}$ normalizes $\bar{G}$ and 
$\genby{G,\Gamma}/L \simeq \genby{\bar{G},\bar{\Gamma}}$. 
Let $M$ be  a $\genby{\bar{G},\bar{\Gamma}}$-module, or equivalently,
a $\genby{G,\Gamma}$-module on which $L$ acts trivially. 
Then $(S_G \otimes M)^G=((S_G)^L \otimes M)^{\bar{G}}$. Hence in view
of (\ref{tau iso})
we have an isomorphism of $\Gamma$-modules 
(on which $\Gamma \cap L$ acts trivially)
\begin{equation}\label{iso quotient}
(\bar{S}_{\bar{G}} \otimes M)^{\bar{G}}
\overset{\bar\tau\otimes\Id_M}{-\!\!\!-\!\!\!-\!\!\!\lr}(S_G \otimes M)^G.
\end{equation}

In (\ref{iso quotient}), we use the $G$-equivariance of $\bar\tau$ to 
restrict to the $G$-fixed points, noting that on the left, $G$ acts via
$\bar G$, since $L$ acts trivially. We denote the map of (\ref{iso quotient})
by $\tau_M$.
Similarly, $\tau_M$ will also denote the isomorphism of $\Gamma$-modules
$(\bar{\CH} \otimes M)^{\bar{G}} \to (\CH \otimes M)^G$ (cf. Corollary \ref{harmonic}).

Next assume that the element $\gamma \in \Gamma$ acts 
semisimply on $M$. Let $\bar{\gamma}$ be its image 
in $\bar{\Gamma}$. Then $\bar{\gamma}$ is semisimple, and the next 
statement follows easily from the above remarks.  

\begin{lemma}\label{compare bases}
Let $\CB(M,\bar{\gamma})$ be a basis 
of $(\bar{\CH} \otimes M)^{\bar{G}}$ consisting of $\bar{\gamma}$-eigenvectors. 
Then $(\tau_M(\iota))_{\iota \in \CB(M,\bar{\gamma})}$ is a basis 
of $(\CH \otimes M)^G$ consisting of $\gamma$-eigenvectors.
\end{lemma}

As an immediate consequence, we have
\begin{corollary}\label{equal factors}
The $M$-factors of $\bar G$ coincide with those of $G$.
\end{corollary}

A further easy consequence of (\ref{compare bases}) is 
\begin{corollary}\label{gutkin quotient}
Let $\Psi_M \in S$ be the polynomial defined  
after Theorem \ref{gutkin}, and let
$\bar{\Psi}_M$ be the element of $\bar{S}$ defined in analogous 
fashion for $M$ as $\bar G$-module. Then 
$\tau(\bar{\Psi}_M) \dotequal \Psi_M$.
\end{corollary}

\subsection{Decomposition of $\bar G$ into graded components}

The vector space $E$ is graded, and $\bar G$ preserves degrees. Therefore
there is a natural decomposition of $\bar G$ into components $\bar G_i$.
In this subsection we relate the various invariants and constants
we have discussed for $\bar G$ to those of the $\bar G_i$. 

Accordingly, write $E=E_1 \oplus E_2 \oplus \cdots $ for the decomposition 
of $E$ into its graded components, with $E_i$ having degree $i$ (so that
$E_i=0$ for all but finitely many $i$). Let $\bar{G}_i$ be the image 
of $\bar{G}$ in $\GL(E_i)$. Since $\bar{G}$ is generated by reflections in $E$, 
the groups $\bar{G}_i$ are generated by reflections in $E_i$ and 
we have $\bar{G} = \bar{G}_1 \times \bar{G}_2 \times \cdots$. 
Let $\bar{S}_{\bar{G}_i}$ be the algebra of coinvariants of $\bar{G}_i$
acting on $E_i$. If $E_i=0$ we take $\bar G_i=1$ and $\bar S_{\bar G_i}=K$.
Then
$$\bar{S}_{\bar{G}} \simeq 
\bar{S}_{\bar{G}_1} \otimes \bar{S}_{\bar{G}_2} \otimes \cdots$$

We begin this subsection with the following observation which relates
eigenvectors for cosets of $G$ to those for $\bar G$ and the $\bar G_i$.  Take
$\gamma \in \Gamma$ and denote by $\bar{\gamma}_i$ its image 
in $\GL(E_i)$, so that 
$\bar{\gamma}=(\bar{\gamma}_1,\bar{\gamma}_2,\cdots)$.

\begin{proposition}\label{eigenvalues} Suppose $\gamma$
has an eigenvector $v$ such that $\gamma v=\zeta v$ for some
$\zeta\in K^\times$. 
Let $\bar{v}$ denote the image of 
$v$ in $V/L$ and write $\bar{v}=\bar v_1 \oplus \bar{v}_2 \oplus \cdots$ with 
$\bar{v}_i \in E_i$ (recall that $\tau$ defines an embedding of $V/L$ into $E$).
Then for each $i$, $\bar{\gamma}_i(\bar{v}_i) = \zeta^i \bar{v}_i.$
\end{proposition}

\begin{proof}
Let $(Q_1,\dots,Q_s)$ be a homogeneous basis of $E^*$ and suppose
that $Q_i$ has degree  $m_i$. Let $(e_1,\dots,e_s)$ be the basis of
$E$ which is dual to 
$(Q_1,\dots,Q_s)$ and let $\pi : V \to V/L \hookrightarrow E$ 
be the natural morphism. Then, by definition, 
$$\bar{v}=\pi(v)=\sum_{i=1}^s Q_i(v) e_i.$$

Then $\bar{\gamma}(\pi(v))=\pi(\gamma(v))=\pi(\zeta v)$. 
So $$\bar{\gamma}(\pi(v))=\sum_{i=1}^s \zeta^{m_i} Q_i(v) e_i,$$
as required.
\end{proof}

The next statement deals with the question of regularity.

\begin{proposition}\label{regular quotient} Assume further
in (\ref{eigenvalues}), that $v$ is $G$-regular.
Then $\bar{v}$ is regular for $\bar{G}$, and {\it a fortiori} $\bar v_i$
is regular for $\bar G_i$, for each $i$. Thus if $\zeta\in K^\times$ is
regular for $\gamma G$, $\zeta^i$ is regular for $\bar\gamma_i \bar G_i$
\end{proposition}

\begin{proof}
For any polynomial $Q\in \bar S=S(E^*)$ and element $w\in V$, we have 
$$
\tau(Q)(w)=Q(\bar w),
$$
where $\bar w$ denotes the image in $E$ of the $L$-orbit of $w$.
Applying (\ref{gutkin quotient}), it follows that for any 
$\bar G$-module $M$, $\bar\Psi_M(\bar v)=\Psi_M(v)\neq 0$, 
since $\Psi_M$ is always a product of linear forms corresponding to
the hyperplanes of $G$, and $v$ is $G$-regular. In particular, this applies 
to the representation $E$, which proves that $\bar v$ is regular. 
The other statements are clear.
\end{proof}

\begin{remark}
The untwisted part of (\ref{regular quotient}) is easily deduced from 
\cite[Theorem 3.12(iii)]{BBR}, while an untwisted analogue of (\ref{eigenvalues})
was stated without proof in \cite[note added in proof]{BBR}.
\end{remark}

We finish by relating the $M$-degrees and constants of $\bar G$ 
and those of the $\bar G_i$.
Let $\bar{\Gamma}_i$ denote the image of $\Gamma$ in $\GL(E_i)$. 
Then $\genby{\bar{G},\bar{\Gamma}}$ is a subgroup of 
$\genby{\bar{G}_1,\bar{\Gamma}_1} \times \genby{\bar{G}_2,\bar{\Gamma}_2} 
\times \cdots$. 
For each $i$, let $\bar{M}_i$ be a $\genby{\bar{G}_i,\bar{\Gamma}_i}$-module and
take the $\genby{\bar{G},\bar{\Gamma}}$-module $M$ to be
$M= \bar{M}_1 \otimes \bar{M}_2 \otimes \cdots$. Then, by (\ref{iso quotient}), 
we have an isomorphism of $\Gamma$-modules
\begin{equation}\label{tensor quotient}
(S_G \otimes M)^G \simeq (\bar{S}_{\bar{G}_1} \otimes \bar{M}_1)^{\bar G_1} \otimes 
(\bar{S}_{\bar{G}_2} \otimes \bar{M}_2)^{\bar{G}_2} \otimes \cdots
\end{equation}

Now define the {\it fake 
$\gamma$-degree} $F_{M,\gamma}(t)$ of $M$ as the polynomial 
$$F_{M,\gamma}(t)=\sum_{i \ge 0} \Trace(\gamma, ((S_G)_i \otimes M^*)^G)
~t^i=\sum_{\iota \in \CB(M,\gamma)}\e_\iota(M,\gamma) t^{m_\iota}.$$

\begin{proposition} We have
$$F_{M,\gamma}(t)=\prod_{i \ge 0} F_{\bar{M}_i,\bar\gamma_i}(t^i).$$
\end{proposition}

\begin{proof}
It is clear from (\ref{tensor quotient}) that we may take 
$\CB(M,\gamma)=\CB(\bar M_1,\bar\gamma_1)\times\CB(\bar M_2,\bar\gamma_2)\cdots$.
If $\beta_i\in \CB(\bar M_i,\bar\gamma_i)$ ($i=1,2,\cdots$) and 
$\iota=\beta_1\otimes\beta_2\otimes\cdots$, then
$$
\e_\iota(M,\gamma) t^{m_\iota}=\prod_i\e_{\beta_i}(\bar M_i,\bar\gamma_i)t^{im_{\beta_i}},
$$
from which the statement is clear.
\end{proof}

This last result is a twisted version of the assertion made without proof
in \cite[note added in proof]{BBR}, which corresponds to the case
$\gamma=\Id$.


\section*{Appendix 1; A list of reflection cosets}

In  this section,  we shall  classify the  reflection cosets  $G\gamma$
where $G$ is irreducible and  $\gamma$ induces a non-inner automorphism
of $G$  (up to multiplication  by scalars) and regular  eigenvalues for
our  choices  of $\gamma$.  The  list  of  reflection cosets  as  above
is  given in  \cite[3.13]{spets}.  The  result is  given  in the  table
concluding this article.  The table in \cite{spets} gives  the image of
$G\gamma$ in the  group of outer automorphisms of  $G$, which describes
the  coset  up  to  a  scalar.  In  each  case  we  choose  a  specific
representative.  First,  for  each  natural number  $d$,  we  choose  a
primitive  $d^{\mathrm{th}}$-root of  unity $\zeta_d$.  We also  assume
that  $\zeta_{de}^e=\zeta_d$  for every  $d$,  $e$.  Before giving  the
table, we explain our conventions and  explain how we get the numerical
results. First, $o(\zeta)$ denotes the order of $\zeta$. Except for the
first two examples in the table, the degrees and codegrees are given in
increasing and decreasing order respectively.

\medskip

\noindent{\bf A product formula.} (cf. \cite{chevie})
The formula
$$\prod_{g\in
G}\det(1-Tg\gamma)=  \prod_i(1-\eps_i   T^{d_i})^{|G|/d_i}$$
is deduced from  the case $M=V$  of (\ref{OS2}) in the  same way
that (1.9) is deduced from (1.8) in \cite{broue}. It is used
in the {\tt CHEVIE} package to compute the $\eps_i$'s.

\medskip

\noindent{\bf The imprimitive groups.}
Let $d$ and $e$ be two natural numbers.
Let $\mu_d$ denote the group of $d$-th  roots of unity in $K$. We choose
a basis $v_i$ of  $V$ such that the group $G(de,e,r)$  is realized as the
group  of  monomial  matrices  with non-zero  entries  in  $\mu_{ed}$ of
which the product of the non-zero  entries lies in $\mu_d$.  The  automorphism
$\gamma$ is induced by the diagonal matrix with  diagonal entries
$(\zeta_{e'd},1,\ldots,1)$,  where $e'$ divides $e$.

If  $\{X_i\}$  is  the basis of $V^*$ dual to  $\{v_i\}$,  the  invariants
of $G(de,e,r)$ are
$P_k=\sum_{j_1<\ldots<j_k}X_{j_1}^{de}\ldots   X_{j_k}^{de}$   for
$k=1,\ldots,r-1$ and $P_r=(X_1\ldots  X_r)^d$. The corresponding degrees
are $ed,2ed,\ldots,(r-1)ed$ and $rd$ and the corresponding $\eps_i$ are
$1,\ldots,1$ and $\zeta_{e'}\inv$.

\smallskip

\noindent{\it The case $d > 1$.}
Let us  determine  the $\eps_i^*$ when $G=G(de,e,r)$ with $d>1$  and
$e>1$.  The codegrees  are
$0,de,\ldots,(r-1)de$. According  to \cite[B.1 ($2'$)]{OT} one  may choose
as  a  basis  of  $(S\otimes V)^G$  the  vectors  $\theta_i=\sum_{j=1}^r
X_j^{(i-1)de+1}\otimes v_j$.  This basis  is $\gamma$-invariant  for our
choice of  $\gamma$, so we  get $\eps_i^*=1$.  We find that  $\zeta$ is
regular when $\zeta^{rd}=\zeta_{e'}$.

\smallskip

\noindent{\it The case $d=1$.}
The group $G(e,e,r)$ is well-generated. Its codegrees
are $0,e,\ldots,(r-2)e$ and  $(r-1)e-r$. We may exclude  the cases $e=1$
where $\gamma$ is inner and $e=r=2$ which is a non-irreducible group, so
we have  $r<(r-1)e$ thus the largest degree is $(r-1)e$. We  may thus
use the relation $\eps_i\eps^*_{\sigma(i)}=\eps_r=1$ to determine the
$\eps_i^*$ which are, ordering the codegrees as above, $1,\ldots,1$ and
$\zeta_{e'}$. We  find that $\zeta$ is  regular whenever $\zeta^{(r-1)e}=1$
or $\zeta^r=\zeta_{e'}$. Note that $\gamma$ is $1$-regular.

\medskip

\noindent{\bf The case ${\boldsymbol{\lexp 3G(4,2,2)}}$.}
Let $\gamma$ be a reflection in $G_6$ of order $3$ and
assume here that $G=G(4,2,2)$ embedded in $G_6$ as a normal
subgroup of index $3$. We have $d_1=d_2=4$ and, since $\gamma$
acts as a reflection on $Y$ (see the proof of Proposition
\ref{ref}), we have $(\eps_1,\eps_2)=(1,\zeta_3)$ or $(1,\zeta_3^2)$.
If we take
$$\gamma = \frac{\zeta_4+1}{2\zeta_3} \begin{pmatrix}
                                -1&1\\
                                \zeta_4&\zeta_4
				\end{pmatrix}.
$$
Then  $\det  \gamma=\zeta_3$,  $\gamma$  stabilises  the  vector  space 
$S_4^G$  which   is  generated  by  $P_1=X_1^4+X_2^4$   and  $P_2=X_1^2 
X_2^2$. An  easy computation shows  that $(\eps_1,\eps_2)=(1,\zeta_3)$. 
Note    that    $\gamma$    permutes    cyclically    the    generators 
$\begin{pmatrix}-1&0\\0&1\end{pmatrix}$,   $\begin{pmatrix}0&-\zeta_4\\ 
\zeta_4&0\end{pmatrix}$ and  $\begin{pmatrix}0&1\\ 1&0\end{pmatrix}$ of 
$G$.                                                                    

On the other hand, $(d_1^*,d_2^*)=(0,4)$ and, since
$\gamma$ is $1$-regular (by direct check),
we find that $\eps_1^*=1$ and $\eps_2^*=\zeta_3^{-1}$.
Now, $\zeta$ is a regular eigenvalue for $G\gamma$ if and only if
$\zeta^4=1$.


\medskip

\noindent{\bf The case ${\boldsymbol{\lexp 2G_7}}$.}
Note that  $\gamma$  comes from  the  normal  embedding
$G_7\subset G_{15}$. We choose $\gamma$ to be a reflection of order $2$.
It acts as a reflection on $Y$. So,
the  pairs $(d_i,\eps_i)$  are  $(12,1),(12,-1)$.
A  direct check shows that the chosen $\gamma$
is $1$-regular, whence the $(d_i^*,\eps_i^*)$ must  be $(0,1),(12,-1)$.
An eigenvalue $\zeta$ is regular if and only if $\zeta^{12}=1$.

\medskip

\noindent{\bf The cases ${\boldsymbol{\lexp 2F_4}}$ and ${\boldsymbol{\lexp 3D_4}}$.}
These  are  Coxeter  groups.  One  may  choose  $\gamma$  as  a  diagram
automorphism  which is  $1$-regular, then  the $\eps_i$  are determined
from $\gamma$'s eigenvalues on $V$. The $\eps_i^*$ are equal to the $\eps_i$
if we order the codegrees in increasing order as well as the degrees.

For $\lexp3D_4$ the pairs $(d_i,\eps_i)$ are
$(2,1),(4,\zeta_3),(4,\zeta_3^2),(6,1)$ where $\zeta_3$ is a primitive cubic
root of unity, and the pairs $(d_i^*,\eps_i^*)$ are
$(0,1),(2,\zeta_3),(2,\zeta_3^2),(4,1)$. An eigenvalue $\zeta$ is regular if
and only if it has order $1,2,3,6$ or $12$.

For $\lexp 2F_4$ the pairs $(d_i,\eps_i)$ are
$(2,1),(6,-1),(8,1),(12,-1)$ and the pairs $(d_i^*,\eps_i^*)$ are
$(0,1),(4,-1),(6,1),(10,-1)$.
An eigenvalue $\zeta$ is regular if and only if it has order $1,2,4,8,12$ or $24$.

\medskip

\noindent{\bf The case ${\boldsymbol{\lexp 4G(3,3,3)}}$.}
In   the  basis   as  above   for   the  imprimitive   groups,  we   may
choose $\gamma=\frac{-1}{\sqrt{-3}}\begin{pmatrix} \zeta_3&1&\zeta_3^2\\
1&1&1\\ \zeta_3^2&1&\zeta_3\end{pmatrix}$.  It is of order  $4$, and does
not stabilize any  set of generators of $G=G(3,3,3)$ of  cardinality $3$.
We find
that the pairs $(d_i,\eps_i)$ are $(3,\zeta_4),(3,-\zeta_4),(6,1)$.
Since the group
is  well-generated we  deduce that the pairs  $(d_i^*,\eps_i^*)$  are
$(0,1),(3,\zeta_4),(3,-\zeta_4)$. An eigenvalue $\zeta$ is regular if
and only if $\zeta^6=1$.

\medskip

\noindent{\bf The case ${\boldsymbol{\lexp 2G(3,3,3)}}$.}
Take for $\gamma$ the square of the above matrix. Then the new $\eps_i$ are
the squares of the previous ones and similarly for the $\eps_i^*$. Again, $
\zeta$ is regular if and only if $\zeta^6=1$.

\medskip

\noindent{\bf The case ${\boldsymbol{\lexp 2G_5}}$.}
Let $\gamma$ be a reflection of order $2$ in $G_{14}$ which does
not lie in $G_5$ (which is a normal subgroup of $G_{14}$ of index $2$).
Again, $\gamma$ acts as a reflection of order $2$ on $Y$. Since
the degrees of $G_5$ are $(6,12)$ and those of $G_{14}$ are $(6,24)$,
we get that $\eps_1=1$ and $\eps_2=-1$.
As  $G$ is well-generated, we may deduce
that the  pairs $(d_i^*,\eps_i^*)$  are  $(0,1),(6,-1)$. We  find that
$\zeta$ is regular if and only if it is of order $1,2,3,6,8$ or $24$.

To obtain the above statement, one may choose as generators
of $G_5$ the elements $s_+$ and $s_-$, where 
$s_\eps=\frac    12\begin{pmatrix}   (-1+\sqrt{-2})\zeta_3&\eps\zeta_{12}\\
\eps\zeta_{12}&(-1-\sqrt{-2})\zeta_3\end{pmatrix}$
for $\eps \in \{+,-\}$. Then, one may take
$\gamma=\begin{pmatrix}1&0\\  0&-1\end{pmatrix}$.
With these choices, $\gamma$ interchanges $s_+$ and $s_-$.

\def\espace{\vphantom{$\begin{array}{c} A \\ A \\ A\end{array}$}}
\def\valeurs#1#2{$\begin{array}{rrrrrr} #1 \\ #2 \end{array}$}
\def\valeursc#1#2{$\begin{array}{cccccc} #1 \\ #2 \end{array}$}
\bigskip

\begin{sidewaystable}
\caption{A table of cosets}
\begin{centerline}{\begin{tabular}{|c|c|c|c|}
\hline
\espace
$(G,\gamma)$ & \valeurs{d_1,\dots,d_r}{\eps_1,\dots,\eps_r} &
\valeurs{d_1^*,\dots,d_r^*}{\eps_1^*,\dots,\eps_r^*}       & $\zeta$ regular \\
\hline
\hline
\espace
\valeursc{\lexp{e'}{G(de,e,r)}}{d > 1}
& \valeurs{ed,&\hspace{-3mm} 2ed,&\hspace{-3mm}\dots,&\hspace{-3mm}(r-1)ed,
&rd}{1,&1,&\hspace{-3mm}\dots,&1,&\hspace{-3mm}\zeta_{e'}^{-1}}
& \valeurs{0,&\hspace{-3mm}ed,&\hspace{-3mm}2ed,&\hspace{-3mm}\dots,
&\hspace{-3mm}(r-1)ed}{1,&1,&1,&\hspace{-3mm}\dots,&1} &
$\zeta^{rd}=\zeta_{e'}$ \\
\hline
\espace
$\lexp{e'}{G(e,e,r)}$
& \valeurs{e,&\hspace{-3mm}2e,&\hspace{-3mm}\ldots,
&\hspace{-3mm}(r-1)e,&r}{1,&1,&\hspace{-3mm}\ldots,&1,&\hspace{-3mm}\zeta_{e'}^{-1}}
& \valeurs{0,&\hspace{-3mm}e,&
\hspace{-3mm}\dots,&\hspace{-3mm}(r-2)e,&
\hspace{-3mm}(r-1)e-r}{1,&1,&\hspace{-3mm}\dots,&1,&\zeta_{e'}} &
\valeursc{\zeta^{rd}=\zeta_{e'}}{\text{or } \zeta^{(r-1)e}=1} \\
\hline
\espace
$\lexp{4}{G(3,3,3)}$ & \valeurs{4, &4, &6}{\zeta_4,&\zeta_4^{-1},&1} &
\valeurs{0,&3,&3}{1,&\zeta_4,&\zeta_4^{-1}} & $\zeta^6=1$ \\
\hline
\espace
$\lexp{2}{G(3,3,3)}$ & \valeurs{4, &4, &6}{-1,&-1,&1} &
\valeurs{0,&3,&3}{1,&-1,&-1} & $\zeta^6=1$ \\
\hline
\espace
$\lexp{3}{G(4,2,2)}$ & \valeurs{4,&4}{\zeta_3,&\zeta_3^2} &
\valeurs{0,&4}{1,&1} & $\zeta^4=1$ \\
\hline
\espace
$\lexp{3}{D_4}$ & \valeurs{2,&4,&4,&6}{1,&\zeta_3,&\zeta_3^2,&1} &
\valeurs{0,&2,&2,&4}{1,&\zeta_3,&\zeta_3^2,&1} &
$o(\zeta)\in \{1,2,3,6,12\}$ \\
\hline
\espace
$\lexp{2}{G_5}$ & \valeurs{6,&12}{1,&-1} & \valeurs{0,&6}{1,&-1} &
$o(\zeta)\in\{1,2,3,6,8,24\}$ \\
\hline
\espace
$\lexp{2}{G_7}$ & \valeurs{12,&12}{1,&-1} & \valeurs{0,&12}{1,&-1} &
$\zeta^{12}=1$ \\
\hline
\espace
$\lexp{2}{F_4}$ & \valeurs{2,&6,&8,&12}{1,&-1,&1,&-1} &
\valeurs{0,&4,&6,&10}{1,&-1,&1,&-1} &
$o(\zeta) \in \{1,2,4,8,12,24\}$ \\
\hline
\end{tabular}}\end{centerline}
\end{sidewaystable}

\section*{Appendix 2; Proofs of (\ref{gutkin}) and (\ref{pregalois})}

\begin{proof}[Proof of Theorem \ref{gutkin}]
Let $\Gut_M=\prod_{H\in\CA}L_H^{N_H(M)}$ and let $\Gut_M'$ be the element
of $S$ defined by $\prod_{\iota \in \CB(M)} \iota =
\Gut_M' \otimes (y_1\wedge \ldots\wedge  y_m)$. Observe that
\begin{equation}\label{gutreduction}
\Gut_{M \oplus M'} = \Gut_M \Gut_{M'} \text{  and  }
\Gut_{M \oplus M'}' = \Gut_M' \Gut_{M'}'.
\end{equation}
Now, let us first assume that $\CA(G)=\{H\}$, so that $G=G_H$.
By (\ref{gutreduction}), the Theorem need only be checked when
$M$ is irreducible, i.e. affords the character $\det^i$ for some
$i \in \{0,1,\dots,e-1\}$, where $e=|G|$. Since here
$\CH$ has basis the set of $L_H^i$, this case is clear.

Now, consider  the general  case. Let $(u_1,\dots,u_m)$  be a  basis of 
$(\CH  \otimes M^*)^G$.  Write $u_i=\sum_{j=1}^m  q_{ji} \otimes  y_j$, 
with  $q_{ji} \in  S$.  Then $u_1  \dots  u_m =  \det  (q_{ij}) .  (y_1 
\wedge \dots  \wedge y_m)$, so  that $\Gut_M'  = \det q_{ij}$.  Now for 
every  $H\in \CA(G)$,  $u_i \in  (S \otimes  M^*)^G \subset  (S \otimes 
M^*)^{G_H}$.  Hence  if  we  express  $u_i$  as  a  linear  combination 
of  elements of  an  $S^{G_H}$-basis of  $(S  \otimes M^*)^{G_H}$  with 
coefficients in  $S^{G_H}$, we see  that $\Gut_M'$ is divisible  by its 
analogue  for $G_H$,  which  $\dotequal L_H^{N_H(M)}$  by the  previous 
discussion.  Since  the distinct  $L_H$  are pairwise coprime, it  follows  that 
$\Gut_M$  divides  $\Gut_M'$.  It   therefore  suffices  to  show  that 
$\Gut_M'$ is non-zero, and has degree $\sum_{H \in \CA(G)} N_H(M)$.     

For the first statement, we prove that if $v\in V$ is such that
$C_G(v)=1$, where $C_G(v)=\{x\in G\mid x(v)=v\}$,
then $\Gut_M'(v)\neq 0$. Note that for any non-trivial element $x\in G$, 
$V_x:=\{v\in V\mid x(v)=v\}$ is a subspace of positive codimension in $V$,
whence $V\setminus \bigcup_{x\in G, x\neq \Id_V} V_x\neq \emptyset$, 
so that such elements $v$ exist. Given one,
if $Gv$ is its $G$-orbit, 
then the map $G \to Gv$ defined by $g \mapsto g(v)$ is bijective. 
Let  $\CF$ be  the space  of  functions $Gv  \to K$,  endowed with  its 
natural  $G$-module  structure. The  restriction  map  $S \to  \CF$  is 
evidently surjective. Since elements of $S^G$ are constant on $Gv$, the 
restriction map $\CH \to \CF$ is  also surjective, and by dimension, is 
an isomorphism  of $G$-modules. Let $f  \in \CH$ be such  that $f(v)=1$ 
and $f(g(v))=0$ if $g \neq \Id_V$.  Consider the map $\nu : \CH \otimes 
M^* \to M^*$, given by $h \otimes x \mapsto h(v) x$. The restriction of 
$\nu:(\CH  \otimes  M^*)^G\to  M^*$  is an  isomorphism  of  $K$-vector 
spaces, since by dimension, it suffices  to show that it is surjective, 
which easily  follows from the fact  that $\nu(\sum_{g \in G}  \lexp gf 
\otimes \lexp gx) =  x$ for every $x \in M^*$. But  the matrix of $\nu$ 
with respect  to the  bases $(u_1,\dots,u_m)$ and  $(y_1,\dots,y_m)$ is 
exactly $q_{ji}(v)$, whose determinant is $\Gut_M'(v)$; it follows that
$\Gut_M'$ is non-zero.                 

It   remains  to   show   that  $N(M)=\sum_{H   \in  \CA(G)}   N_H(M)$. 
Define  the  {\it fake  degree}  $F_M(t)$  of  $M$ as  the  Poincar\'e 
polynomial  $\sum_{\iota \in  \CB(M)} t^{m_\iota}$  of $(\CH\otimes 
M^*)^G$,  where  $t$  is  an  indeterminate.  It  is  then  clear  that 
$N(M)=\frac{\partial   F_M}{\partial   t}\mid_{t=1}$.   However   from 
Molien's  formula,  if  $\chi$  is  the  character  of  the  $G$-module 
$M$, we have  $F_M(t)=\frac{\prod_{\iota \in \CB(V)}(1-t^{d_\iota})}{|G|} 
\sum_{g\in G}\frac{\chi(g)}{\det_V(1-gt)}$. Taking  derivatives, we see 
that the terms where  $g$ is not $1$ or a  reflection do not contribute 
to $N(M)$. To sum the remaining terms,  we use the fact that the number 
$|\Ref(G)|$ of reflections of $G$ is equal to $\sum_\iota{(d_\iota-1)}$ 
to obtain                                                               
$N(M)=\chi(1)|\Ref(G)|/2+\sum_{s\in\Ref(G)}\frac{\chi(s)}{\det(s\mid
V)-1}$.  This  expression  is  exactly the  sum  of  the  corresponding 
expressions  for $G_H$,  over  all $H\in  \CA(G)$,  whence the  result. 
\end{proof}                                                             

\begin{proof}[Proof of Lemma \ref{pregalois}]
Let $M$  be a  $G$-module of  dimension $m$ such  that any  reflection of 
$G$  acts as  a  reflection in  $M$.  Then for  $H\in\CA$, we have 
in the notation preceding the statement of Theorem \ref{gutkin}, the
$G_H$-module decomposition $M^*=\oplus_{i=1}^m\xi_H^{e_i}$,
where $e_1\neq 0$ and $e_i=0$ for $i>1$. Then $\Lambda^m  M^*=\xi_H^{e_1}$,
and $N_H(M)=N_H(\Lambda^m  M)=e_1$. Since this holds for any
$H\in \CA$, and $N(M)=\sum_{H\in\CA}N_H(M)$ for any $M$, the Lemma follows.                                            
\end{proof}

\end{document}